# Strict Fréchet Differentiability of the Metric Projection Operator in Hilbert Spaces

Jinlu Li

Department of Mathematics
Shawnee State University
Portsmouth, Ohio 45662
USA

**Abstract**

In this paper, we prove strict Fréchet differentiability of the metric projection operator onto closed balls in Hilbert spaces and onto positive cones in Euclidean spaces. We find the exact expressions for Fréchet derivatives. Since Fréchet differentiability implies Gâteaux directional differentiability, the results obtained in this paper strengthen the results obtained in [8] and [10] about the directional differentiability of the metric projection operator onto closed balls in Hilbert spaces and positive cones in Euclidean spaces and in the real Hilbert space $l_2$.

1. **Introduction**

Throughout this paper, let $(H, \|\cdot\|)$ be a real Hilbert space with inner product $\langle \cdot, \cdot \rangle$ and with the origin $\theta$. Let $C$ be a nonempty closed and convex subset of $H$. Let $P_C: H \to C$ denote the (standard) metric projection operator satisfying

$$\|x - P_C(x)\| \leq \|x - z\|, \text{ for all } z \in C.$$

It is well-known that $P_C$ is a well-defined single-valued mapping. In operator theory in Hilbert spaces, one of the most important and most useful operators is the metric projection operator. It has many useful properties, such as continuity, monotonicity, and non-expensiveness. These properties have been studied by many authors and have been applied to many fields in analysis in Hilbert spaces (see 1, 5, 9, 17).

In addition to continuity, the smoothness of mappings in Hilbert spaces has been studied with respect to several types of differentiability (see [2−4, 6−7, 10, 14−16]), which have been applied to approximation theory, optimization theory, and variational inequalities in Hilbert spaces. For details, see [11−13].

In particular, the Gâteaux directional differentiability of the metric projection operator $P_C$, at a given point and along a certain direction, has been introduced; this has been defined in uniformly convex and uniformly smooth Banach spaces in [8] (we recall this definition in section 2). In [8] and [10], some properties of Gâteaux directional derivatives of $P_C$ in both Banach spaces and Hilbert spaces have been proved. In particular, when $C$ is a closed ball, or a closed and convex cone, or a closed and convex cylinder, the exact Gâteaux directional derivatives of $P_C$ are provided in [8] and [10].

Compared to the Gâteaux directional differentiability of $P_C$, a stricter differentiability of $P_C$ is Fréchet differentiability. A more restrictive concept of Fréchet differentiability is strict Fréchet differentiability of the metric projection $P_C$. The two concepts of Fréchet differentiability will be recalled in section 2.

In general, let $f: H \to H$ be a single-valued mapping. In Proposition 2.1 in section 2 of this paper, we prove that all these three concepts of differentiability of a single-valued mapping $f$ satisfy that, for any given point $\bar{x} \in H$, one has

$f$ is strictly Fréchet differentiable at $\bar{x}$

$\Longrightarrow f$ is Fréchet differentiable at $\bar{x}$

$\Longrightarrow f$ is Gâteaux directionally differentiable at $\bar{x}$.

However, the converse statements may not hold, which will be demonstrated by the metric projection operator as follows. Both part (iii) of Theorem 3.3 in section 3, part (iv) of Theorem 4.4 in section 4, Theorem 5.3 of this paper independently show that, for some point $\bar{x}$, we have

The metric projection is Gâteaux directionally differentiable at $\bar{x}$

$\nRightarrow$ The metric projection is Fréchet differentiable at this point $\bar{x}$.

In this paper, we focus on the differentiability of the metric projection operator onto closed balls in Hilbert spaces, the positive cones in Euclidean spaces and in the real Hilbert space $l_2$. By adapting the concepts of Fréchet differentiability and strict Fréchet differentiability of the metric projection used in [13], we strengthen the results from Gâteaux directional differentiability of the metric projection operator obtained in [8] and [10] to strict Fréchet differentiability and Fréchet differentiability of the metric projection operator.

In [13], Mordukhovich introduced the concept of generalized differentiation and provide many useful properties, such as differential calculus in Banach spaces. The theory of generalized differentiation of mappings in Banach spaces has been widely applied to nonlinear analysis presented in this book. Furthermore, some applications of the Gâteaux directional differentiability of $P_C$ has been applied to approximation problems, convex programming, optimal control, and so forth (see [3, 7, 14, 15]).

In section 2, we recall some concepts and properties of differentiability of mappings in Hilbert spaces. In section 3 (Theorem 3.3), we prove the strict Fréchet differentiability of the metric projection operator onto closed balls in Hilbert spaces. We also find the precise solutions of the Fréchet derivatives of the projection operator, where it is strictly Fréchet differentiable. Following Theorem 3.3, we provide some examples to demonstrate this theorem. In section 4, we prove the strict Fréchet differentiability of the metric projection operator onto positive cones in Euclidean spaces and we find the exact representations of the Fréchet derivatives. In section 5, we consider the real Hilbert space $l_2$. In this section, we study the Fréchet differentiability of the metric projection operator onto the positive cone in $l_2$, which is very different from the Fréchet differentiability in Euclidean spaces.

## 2. Preliminaries

### 2.1. The metric projection operator in Hilbert spaces

Since every Hilbert space $H$ is a uniformly convex and uniformly smooth Banach space, the norm $\|\cdot\|$ of $H$ has the following properties.

$$\lim_{t \downarrow 0} \frac{\|x+tv\| - \|x\|}{t} = \frac{\langle x,v \rangle}{\|x\|}, \text{ for any } (x, v) \in H\setminus\{\theta\} \times H\setminus\{\theta\}.$$

In particular,

$$\lim_{t \downarrow 0} \frac{\|x+tv\| - \|x\|}{t} = \langle x, v \rangle, \text{ uniformly for } (x, v) \in \mathbb{S} \times \mathbb{S}. \qquad (1.3)$$

For any $y \in C$, the inverse image of $y$ by the metric projection $P_C$ in $H$ is defined by

$$P_C^{-1}(y) = \{x \in H : P_C(x) = y\}.$$

We provide some properties of the projection operator below listed as a proposition without proof. See [1, 5, 11] for more details and proofs.

**Proposition 2.1**. *Let C be a nonempty closed and convex subset of a Hilbert space H. The metric projection $P_C: X \to C$ has the following properties*

(i)     *$C$ is a Chebyshev set in $H$. That is, for any $x \in H$, there is a unique point $P_C x \in C$;*

(ii)    *The basic variational principle: for any $x \in H$ and $u \in C$,*

$$u = P_C x \iff \langle x - u, u - z \rangle \geq 0, \text{ for all } z \in C;$$

(iii)    *For any $x \in H$ and $u \in C$,*

$$u = P_C x \iff \langle x - u, x - z \rangle \geq \|x - u\|^2, \text{ for all } z \in C;$$

(iv)    *$P_C$ is (strongly) monotone, that is*

$$\langle P_C x - P_C y, x - y \rangle \geq \|P_C x - P_C y\|^2, \text{ for all } x, y \in H;$$

(v)     *Non-expansiveness*:

$$\|P_C x - P_C y\| \leq \|x - y\|, \text{ for any } x, y \in H;$$

(vi)    *Let $c \in H$ and $r > 0$. Then we have*

$$P_{\mathbb{B}(c,r)}(x) = c + \frac{r}{\|x-c\|}(x - c), \text{ for any } x \in H\setminus\mathbb{B}(c, r).$$

       *In particular, we have*

$$P_{\mathbb{B}}(x) = \frac{x}{\|x\|}, \text{ for any } x \in H \text{ with } \|x\| > 1.$$

### 2.2. Differentiability of the metric projection operator onto closed and convex subsets

We first recall the concepts of the differentiability of single-valued mappings from $H$ to $H$.

**Definition 3.1 in [8].** (Gâteaux directional differentiability of single-valued mappings in Hilbert spaces). Let $f: H \to H$ be a single-valued mapping. For $x \in H$ and $w \in H$ with $w \neq \theta$, if the following limit exists, which is a point in $H$,

$$f'(x)(w) = \lim_{t \downarrow 0} \frac{f(x+tw)-f(x)}{t}, \tag{2.1}$$

then, $f$ is said to be (Gâteaux) directionally differentiable at point $x$ along the direction $w$. $f'(x)(w)$ is called the (Gâteaux) directional derivative of $f$ at point $x$ along direction $w$. Let $A$ be a subset in $H$. If $f$ is (Gâteaux) directionally differentiable at every point $x \in A$, then $f$ is said to be (Gâteaux) directionally differentiable on $A \subseteq H$.

**Definition 1.13 in [17].** (Fréchet differentiability of single-valued mappings in Hilbert spaces). Let $f: H \to H$ be a single-valued mapping. For any given $\bar{x} \in H$, if there is a linear continuous mapping $\nabla f(\bar{x}): H \to H$ such that

(i)
$$\lim_{x \to \bar{x}} \frac{f(x)-f(\bar{x})-\nabla f(\bar{x})(x-\bar{x})}{\|x-\bar{x}\|} = \theta. \tag{2.2}$$

then $f$ is said to be Fréchet differentiable at $\bar{x}$ and $\nabla f(\bar{x})$ is called the Fréchet derivative of $f$ at $\bar{x}$;

(ii) More strictly, if

$$\lim_{(u,v) \to (\bar{x},\bar{x})} \frac{f(u)-f(v)-\nabla f(\bar{x})(u-v)}{\|u-v\|} = \lim_{u \to \bar{x}, v \to \bar{x}} \frac{f(u)-f(v)-\nabla f(\bar{x})(u-v)}{\|u-v\|} = \theta, \tag{2.3}$$

then $f$ is said to be strictly Fréchet differentiable at $\bar{x}$.

**Proposition 2.1.** *Let $H$ be a Hilbert space and let $f: H \to H$ be a single-valued mapping. Then, for any given point $\bar{x} \in H$, one has*

(a) *$f$ is strictly Fréchet differentiable at $\bar{x}$*
*$\Longrightarrow$ $f$ is Fréchet differentiable at $\bar{x}$*
*$\Longrightarrow$ $f$ is Gâteaux directionally differentiable at $\bar{x}$ along every direction satisfying*

$$f'(\bar{x})(w) = \nabla f(\bar{x})(w), \text{ for any } w \in H \text{ with } w \neq \theta. \tag{2.4}$$

(b) *$f$ is Gâteaux directionally differentiable at $\bar{x}$ along every direction*
*$\not\Rightarrow$ $f$ is Fréchet differentiable at the given point $\bar{x}$.*

*Proof.* Proof of (a). For any given $\bar{x} \in H$, suppose that $f$ is Fréchet differentiable at $\bar{x}$ with $\nabla f(\bar{x})$ being the Fréchet derivative of $f$ at $\bar{x}$, which is a linear continuous mapping $\nabla f(\bar{x}): H \to H$ such that

$$\lim_{x \to \bar{x}} \frac{f(x)-f(\bar{x})-\nabla f(\bar{x})(x-\bar{x})}{\|x-\bar{x}\|} = \theta.$$

Then, for any given $w \in H$ with $w \neq \theta$, in the limit (2.2), we take a special direction approaching to $\bar{x}$: $\bar{x} + tw$, for $t \downarrow 0$. By the linearity of $\nabla f(\bar{x})$, we have

$$\theta = \lim_{\bar{x}+tw \to \bar{x}} \frac{f(\bar{x}+tw)-f(\bar{x})-\nabla f(\bar{x})((\bar{x}+tw)-\bar{x})}{\|(\bar{x}+tw)-\bar{x}\|}$$

$$= \lim_{t\downarrow 0} \frac{f(\bar{x}+tw)-f(\bar{x})-t\nabla f(\bar{x})(w)}{t\|w\|}$$

$$= \frac{1}{\|w\|}\lim_{t\downarrow 0} \frac{f(\bar{x}+tw)-f(\bar{x})}{t} - \frac{\nabla f(\bar{x})(w)}{\|w\|}.$$

This implies

$$\nabla f(\bar{x})(w) = \lim_{t\downarrow 0} \frac{f(\bar{x}+tw)-f(\bar{x})}{t} = f'(\bar{x})(w).$$

Proof of (b). For any given $r > 0$, let $P_{r\mathbb{B}}: H \to r\mathbb{B}$ be the metric projection operator, which is a well-defined single-valued mapping from $H$ onto $r\mathbb{B}$. Then part (b) could be shown by taking this special mapping $f = P_{r\mathbb{B}}$. From the results of part (iii) of Theorem 5.2 in [8], for any given point $\bar{x} \in r\mathbb{S}$, we have that $P'_{r\mathbb{B}}(\bar{x})(v)$ exists, for any $v \in H\setminus\{\theta\}$. However, from the results in part (iii) of Theorem 3.3 in section 3 of this paper, we see that $\nabla P_{r\mathbb{B}}(\bar{x})$ does not exist, for any point $\bar{x} \in r\mathbb{S}$. □

Let $C$ be a nonempty closed and convex subset of $H$. We consider some properties of the differentiability of $P_C$, which will be used in the proof of Theorems 3.3 and 4.3.

**Proposition 2.2.** *Let $C$ be a nonempty closed and convex subset of $H$. Let $y \in C$. Suppose $(P_C^{-1}(y))^o \neq \emptyset$. Then, $P_C$ is strictly Fréchet differentiable on $(P_C^{-1}(y))^o$ such that,*

$$\nabla P_C(\bar{x}) = \theta, \text{ for any } \bar{x} \in (P_C^{-1}(y))^o.$$

*That is, for any $\bar{x} \in (P_C^{-1}(y))^o$, one has*

$$\nabla P_C(\bar{x})(w) = \theta, \text{ for every } w \in H \text{ with } w \neq \theta.$$

*Proof.* For any $\bar{x} \in (P_C^{-1}(y))^o$, there is a positive number $p$ such that

$$\mathbb{B}(\bar{x},p) \subseteq (P_C^{-1}(y))^o.$$

This implies that

$$\lim_{u\to\bar{x},v\to\bar{x}} \frac{P_C(u)-P_C(v)-\theta}{\|u-v\|}$$

$$= \lim_{\substack{u\to\bar{x},v\to\bar{x} \\ u,v \in \mathbb{B}(\bar{x},p)}} \frac{P_C(u)-P_C(v)}{\|u-v\|}$$

$$= \lim_{\substack{u\to\bar{x},v\to\bar{x} \\ u,v \in \mathbb{B}(\bar{x},p)}} \frac{y-y}{\|u-v\|} = \theta. \qquad \square$$

**Proposition 2.3.** *Let $C$ be a nonempty closed and convex subset of $H$. Suppose $C^o \neq \emptyset$. Then $P_C$ is strictly Fréchet differentiable on $C^o$ satisfying*

$$\nabla P_C(\bar{x}) = I_H, \text{ for any } \bar{x} \in C^o.$$

*That is, for any $\bar{x} \in C^o$, one has*

$$\nabla P_C(\bar{x})(w) = w, \text{ for every } w \in H \text{ with } w \neq \theta.$$

*Proof.* For any $\bar{x} \in C^o$, there is a positive number $p$ such that

$$\mathbb{B}(\bar{x}, p) \subseteq C^o.$$

This implies that

$$\lim_{u \to \bar{x}, v \to \bar{x}} \frac{P_C(u) - P_C(v) - I_H(u-v)}{\|u-v\|}$$

$$= \lim_{\substack{u \to \bar{x}, v \to \bar{x} \\ u,v \in \mathbb{B}(\bar{x},p)}} \frac{P_C(u) - P_C(v) - (u-v)}{\|u-v\|}$$

$$= \lim_{\substack{u \to \bar{x}, v \to \bar{x} \\ u,v \in \mathbb{B}(\bar{x},p)}} \frac{u - v - (u-v)}{\|u-v\|} = \theta. \qquad \square$$

### 3. Strict Fréchet differentiability of the metric projection operator onto balls

Let $\mathbb{B}$ denotes the unit closed ball in a Hilbert space $H$. For any $r > 0$, $r\mathbb{B}$ denotes the closed ball with radius $r$ and centered at $\theta$. Let $\mathbb{S}$ be the unit sphere in $H$. Then, $r\mathbb{S}$ is the sphere in $H$ with radius $r$ and centered at $\theta$. For any $c \in H$ and $r > 0$, let $\mathbb{B}(c, r)$ denote the closed ball in $H$ with radius $r$ and centered at $c$. In this notation, $\mathbb{B}(\theta, 1) = \mathbb{B}$ and $\mathbb{B}(\theta, r) = r\mathbb{B}$. $\mathbb{S}(c, r)$ denotes the sphere in $H$ with center $c$ and with radius $r$. Let $I_H$ denote the identity mapping in $H$.

### 3.1 Review Gâteaux directional differentiability of the metric projection operator

In [8, 10], Gâteaux directional differentiability of $P_{\mathbb{B}(c,r)}$ was studied in both Banach spaces and Hilbert spaces. Before we state the results obtained in [10], we review some related notations.

For any $x \in \mathbb{S}(c, r)$, two subsets $x^{\uparrow}_{(c,r)}$ and $x^{\downarrow}_{(c,r)}$ of $H \setminus \{\theta\}$ are defined by: for $v \in H$ with $v \neq \theta$,

(a) $v \in x^{\uparrow}_{(c,r)} \iff$ there is $\delta > 0$ such that $\|(x + tv) - c\| \geq r$, for all $t \in (0, \delta)$;
(b) $v \in x^{\downarrow}_{(c,r)} \iff$ there is $\delta > 0$ such that $\|(x + tv) - c\| < r$, for all $t \in (0, \delta)$.

In particular, we write

(c) $x^{\uparrow}_r = x^{\uparrow}_{(\theta,r)}$: $v \in x^{\uparrow}_r \iff$ there is $\delta > 0$ such that $\|x + tv\| \geq r$, for all $t \in (0, \delta)$;
(d) $x^{\downarrow}_r = x^{\downarrow}_{(\theta,r)}$: $v \in x^{\downarrow}_r \iff$ there is $\delta > 0$ such that $\|x + tv\| < r$, for all $t \in (0, \delta)$.

**Theorem 4.2 in [10] or Theorem 5.2 in [8].** *Let $C = \mathbb{B}(c, r)$ be a closed ball in Hilbert space $H$. Then, $P_C$ is directionally differentiable on $H$ such that, for every $w \in H$ with $w \neq \theta$, $P_C$ has the following Gâteaux directional differentiability properties.*

(i)      *For any $x \in \mathbb{B}(c,r)^o$, we have*

    (a)      $P'_C(x)(w) = w,$
    (b)      $P'_C(x)(x) = x,\ \text{for } x \neq \theta;$

(ii)    *For any $x \in H\setminus\mathbb{B}(c,r)$, we have*

    (a)      $P'_C(x)(w) = \frac{r}{\|x-c\|^3}(\|x-c\|^2 w - \langle x-c, w\rangle(x-c)),$
    (b)      $P'_C(x)(x-c) = \theta;$

(iii)    *For any $x \in \mathbb{S}(c,r)$, we have*

    (a)      $P'_C(x)(w) = w - \frac{1}{r^2}\langle x-c, w\rangle(x-c),\ \text{if } w \in x^\uparrow_{(c,r)};$
    (b)      $P'_C(x)(x-c) = \theta;$
    (c)      $P'_C(x)(w) = w,\ \text{if } w \in x^\downarrow_{(c,r)}.$

## 3.2 Strict Fréchet differentiability of metric projection onto balls in Hilbert spaces

For any $x, \bar{x} \in H$, as usual, we write $x \perp \bar{x}$ if and only if $\langle x, \bar{x}\rangle = 0$. For any $\bar{x} \in H\setminus\{\theta\}$, let $S(\bar{x})$ denote the one-dimensional subspace of $H$ generated by $\bar{x}$. Let $O(\bar{x})$ denote the orthogonal subspace of $\bar{x}$ (or $S(\bar{x})$) in $H$. $H$ has the following orthogonal decomposition

$$H = S(\bar{x}) \oplus O(\bar{x}).$$

More precisely speaking, for this given $\bar{x} \in H\setminus\{\theta\}$ and for any $x \in H$, $x$ enjoys the following orthogonal representation

$$x = \frac{\langle x,\bar{x}\rangle}{\|\bar{x}\|^2}\bar{x} + \left(x - \frac{\langle x,\bar{x}\rangle}{\|\bar{x}\|^2}\bar{x}\right),\ \text{for all } x \in H. \tag{3.1}$$

By using the above orthogonal representations of elements in $H$, for this given fixed $\bar{x} \in H\setminus\{\theta\}$, we define a real valued function $a(\bar{x};\cdot)\colon H \to \mathbb{R}$ by

$$a(\bar{x}; x) := \frac{\langle x,\bar{x}\rangle}{\|\bar{x}\|^2},\ \text{for all } x \in H.$$

And a mapping $o(\bar{x};\cdot)\colon H \to O(\bar{x})$ by

$$o(\bar{x}; x) := x - \frac{\langle x,\bar{x}\rangle}{\|\bar{x}\|^2}\bar{x},\ \text{for all } x \in H.$$

The following lemma provides some properties of $a(\bar{x};\cdot)$ and $o(\bar{x};\cdot)$. These properties will play important roles and will be repeatedly used in the proof of Theorem 3.3.

**Lemma 3.1**. *For any given fixed $\bar{x} \in H\setminus\{\theta\}$, the real valued function $a(\bar{x};\cdot)$ and the mapping $o(\bar{x};\cdot)$ have the following properties.*

    (i)      *$a(\bar{x};\cdot)\colon H \to \mathbb{R}$ is a real valued linear and continuous function;*

(ii)  $o(\bar{x}; \cdot): H \to O(\bar{x})$ is a linear and continuous mapping;
(iii) $a(\bar{x}; \cdot)\bar{x}$ and $o(\bar{x}; \cdot)$ are orthogonal of each other and for any $u, v \in H$, we have

(a) $\langle u, v \rangle = a(\bar{x}; u)a(\bar{x}; v)\|\bar{x}\|^2 + \langle o(\bar{x}; u), o(\bar{x}; v) \rangle$;
(b) $\langle a(\bar{x}; u)\bar{x}, o(\bar{x}; v) \rangle = 0$;
(c) $\|a(\bar{x}; u)\bar{x} + o(\bar{x}; v)\|^2 = (a(\bar{x}; u))^2\|\bar{x}\|^2 + \|o(\bar{x}; v)\|^2$;
(d) $\|u\|^2 = (a(\bar{x}; u))^2\|\bar{x}\|^2 + \|o(\bar{x}; u)\|^2$;
(e) $\|u + v\|^2 = (a(\bar{x}; u) + a(\bar{x}; u))^2\|\bar{x}\|^2 + \|o(\bar{x}; u) + o(\bar{x}; v)\|^2$.

(iv) $\qquad u \to \bar{x} \iff a(\bar{x}; u) \to 1 \text{ and } o(\bar{x}; u) \to \theta, \text{ for } u \in H.$ $\qquad$ (3.2)

*Proof.* This lemma can be straightforwardly checked and the proof is omitted here. □

**Remarks 3.2.** Both $a(\bar{x}; \cdot)$ and $o(\bar{x}; \cdot)$ depend on $\bar{x}$. However, since the proof of Theorem 3.3 is considered to be long and both $a(\bar{x}; \cdot)$ and $o(\bar{x}; \cdot)$ are repeatedly used in the proof, for the sake of simplicity, $a(\bar{x}; \cdot)$ and $o(\bar{x}; \cdot)$ are abbreviated as $a(\cdot)$ and $o(\cdot)$, respectively.

The following theorem is to prove the strict Fréchet differentiability of the metric projection onto closed balls in Hilbert spaces.

**Theorem 3.3.** *Let $H$ be a Hilbert space. For any $r > 0$, the metric projection $P_{r\mathbb{B}}: H \to r\mathbb{B}$ has the following differentiability properties.*

(i) $\quad P_{r\mathbb{B}}$ *is strictly Fréchet differentiable on $r\mathbb{B}^o$ satisfying*

$$\nabla P_{r\mathbb{B}}(\bar{x}) = I_H, \text{ for every } \bar{x} \in r\mathbb{B}^o.$$

*That is,*

$$\bar{x} \in r\mathbb{B}^o \implies \nabla P_{r\mathbb{B}}(\bar{x})(x) = x, \text{ for every } x \in H.$$

(ii) $\quad P_{r\mathbb{B}}$ *is strictly Fréchet differentiable on $H \backslash r\mathbb{B}$ such that, for every $\bar{x} \in H \backslash r\mathbb{B}$,*

$$\nabla P_{r\mathbb{B}}(\bar{x})(x) = \frac{r}{\|\bar{x}\|}\left(x - \frac{\langle x, \bar{x} \rangle}{\|\bar{x}\|^2}\bar{x}\right), \text{ for every } x \in H.$$

*In particular, we have*

(a) $\quad \nabla P_{r\mathbb{B}}(\bar{x})(x) = \frac{r}{\|\bar{x}\|}x, \text{ if } x \perp \bar{x}, \text{ for } x \in H$;
(b) $\quad \nabla P_{r\mathbb{B}}(\bar{x})(\bar{x}) = \theta.$

(iii) *On the subset $r\mathbb{S}$, we have*

(I) $P_{r\mathbb{B}}$ *is Gâteaux directionally differentiable on $r\mathbb{S}$ satisfying that, for every point $\bar{x} \in r\mathbb{S}$, the following representations are satisfied*

(a) $\quad P'_{r\mathbb{B}}(\bar{x})(w) = w - \frac{1}{r^2}\langle \bar{x}, w \rangle x, \text{ if } w \in \bar{x}^\uparrow_r;$

  (b) $P'_{r\mathbb{B}}(\bar{x})(\bar{x}) = \theta$;
  (c) $P'_{r\mathbb{B}}(\bar{x})(w) = w$, if $w \in \bar{x}_r^{\downarrow}$.

 (II) $P_{r\mathbb{B}}$ is not Fréchet differentiable at any point $\bar{x} \in r\mathbb{S}$. That is,

$$\nabla P_{r\mathbb{B}}(\bar{x}) \text{ does not exist, for any } \bar{x} \in r\mathbb{S}.$$

*Proof.* Proof of (i). By (2.3), for any given $\bar{x} \in r\mathbb{B}^o$, we calculate

$$\lim_{u \to \bar{x}, v \to \bar{x}} \frac{P_{r\mathbb{B}}(u) - P_{r\mathbb{B}}(v) - I_H(u-v)}{\|u-v\|}$$

$$= \lim_{\substack{u \to \bar{x}, v \to \bar{x} \\ u,v \in \mathbb{B}(\bar{x},p)}} \frac{P_{r\mathbb{B}}(u) - P_{r\mathbb{B}}(v) - (u-v)}{\|u-v\|}$$

$$= \lim_{\substack{u \to \bar{x}, v \to \bar{x} \\ u,v \in \mathbb{B}(\bar{x},p)}} \frac{u - v - (u-v)}{\|u-v\|}$$

$$= \theta.$$

Hence, $P_{r\mathbb{B}}$ is strictly Fréchet differentiable at $\bar{x}$, with $\nabla P_{r\mathbb{B}}(\bar{x}) = I_H$, for $\bar{x} \in r\mathbb{B}^o$.

Proof of (ii). Let $\bar{x} \in H \backslash r\mathbb{B}$ be arbitrarily given and fixed with $\|\bar{x}\| > r$. By the definition of $o(\cdot)$, we actually have

$$\nabla P_{r\mathbb{B}}(\bar{x})(x) = \frac{r}{\|\bar{x}\|}\left(x - \frac{\langle x,\bar{x}\rangle}{\|\bar{x}\|^2}\bar{x}\right) = \frac{r}{\|\bar{x}\|}o(x), \text{ for every } x \in H.$$

To prove part (ii) of this theorem, we only need to verify that the above formula satisfies the following equations.

$$\theta = \lim_{u \to \bar{x}, v \to \bar{x}} \frac{P_{r\mathbb{B}}(u) - P_{r\mathbb{B}}(v) - \frac{r}{\|\bar{x}\|}o(u-v)}{\|u-v\|}$$

$$= \lim_{u \to \bar{x}, v \to \bar{x}} \frac{P_{r\mathbb{B}}(a(u)\bar{x} + o(u)) - P_{r\mathbb{B}}(a(v)\bar{x} + o(v)) - \frac{r}{\|\bar{x}\|}o(u-v)}{\|u-v\|}$$

$$= \lim_{u \to \bar{x}, v \to \bar{x}} \frac{P_{r\mathbb{B}}(a(u)\bar{x} + o(u)) - P_{r\mathbb{B}}(a(v)\bar{x} + o(v)) - \frac{r}{\|\bar{x}\|}o(u) + \frac{r}{\|\bar{x}\|}o(v)}{\|u-v\|}$$

Since $\bar{x} \notin r\mathbb{B}$ with $\|\bar{x}\| > r$, there is $p > 0$ such that

$$\mathbb{B}(\bar{x}, p) \cap r\mathbb{B} = \emptyset.$$

This implies

$$\lim_{u \to \bar{x}, v \to \bar{x}} \frac{P_{r\mathbb{B}}(a(u)\bar{x} + o(u)) - P_{r\mathbb{B}}(a(v)\bar{x} + o(v)) - \frac{r}{\|\bar{x}\|}o(u-v)}{\|u-v\|}$$

$$= \lim_{\substack{u\to\bar{x},v\to\bar{x} \\ u,v\in\mathbb{B}(\bar{x},p)}} \frac{P_{r\mathbb{B}}(a(u)\bar{x}+o(u))-P_{r\mathbb{B}}(a(v)\bar{x}+o(v))-\frac{r}{\|\bar{x}\|}o(u)+\frac{r}{\|\bar{x}\|}o(v)}{\|u-v\|}$$

$$= \lim_{\substack{u\to\bar{x},v\to\bar{x} \\ u,v\in\mathbb{B}(\bar{x},p)}} \frac{\frac{r}{\|a(u)\bar{x}+o(u)\|}(a(u)\bar{x}+o(u))-\frac{r}{\|a(v)\bar{x}+o(v)\|}(a(v)\bar{x}+o(v))-\frac{r}{\|\bar{x}\|}o(u)+\frac{r}{\|\bar{x}\|}o(v)}{\|u-v\|}$$

$$= r \lim_{\substack{u\to\bar{x},v\to\bar{x} \\ u,v\in\mathbb{B}(\bar{x},p)}} \frac{\frac{a(u)\bar{x}+o(u)}{\|a(u)\bar{x}+o(u)\|}-\frac{a(v)\bar{x}+o(v)}{\|a(v)\bar{x}+o(v)\|}-\frac{o(u)}{\|\bar{x}\|}+\frac{o(v)}{\|\bar{x}\|}}{\|u-v\|}$$

$$= r \lim_{\substack{u\to\bar{x},v\to\bar{x} \\ u,v\in\mathbb{B}(\bar{x},p)}} \frac{\frac{a(u)\bar{x}}{\|a(u)\bar{x}+o(u)\|}+\frac{o(u)}{\|a(u)\bar{x}+o(u)\|}-\frac{o(u)}{\|\bar{x}\|}-\frac{a(v)\bar{x}}{\|a(v)\bar{x}+o(v)\|}-\frac{o(v)}{\|a(v)\bar{x}+o(v)\|}+\frac{o(v)}{\|\bar{x}\|}}{\|u-v\|}$$

$$= r \lim_{\substack{u\to\bar{x},v\to\bar{x} \\ u,v\in\mathbb{B}(\bar{x},p)}} \frac{\frac{a(u)\bar{x}}{\|a(u)\bar{x}+o(u)\|}-\frac{a(v)\bar{x}}{\|a(v)\bar{x}+o(v)\|}+\frac{o(u)}{\|a(u)\bar{x}+o(u)\|}-\frac{o(u)}{\|\bar{x}\|}-\frac{o(v)}{\|a(v)\bar{x}+o(v)\|}+\frac{o(v)}{\|\bar{x}\|}}{\|u-v\|}$$

$$= r \lim_{\substack{u\to\bar{x},v\to\bar{x} \\ u,v\in\mathbb{B}(\bar{x},p)}} \left( \frac{\frac{a(u)\bar{x}}{\|a(u)\bar{x}+o(u)\|}-\frac{a(v)\bar{x}}{\|a(v)\bar{x}+o(v)\|}}{\|u-v\|} + \frac{\left(\frac{o(u)}{\|a(u)\bar{x}+o(u)\|}-\frac{o(u)}{\|\bar{x}\|}\right)-\left(\frac{o(v)}{\|a(v)\bar{x}+o(v)\|}-\frac{o(v)}{\|\bar{x}\|}\right)}{\|u-v\|} \right). \tag{3.3}$$

By (3.2) in Lemma 3.1, we have

$$\begin{aligned} & u\to\bar{x} \quad\Longleftrightarrow\quad a(u)\to 1 \text{ and } o(u)\to\theta. \\ \text{and} \quad & v\to\bar{x} \quad\Longleftrightarrow\quad a(v)\to 1 \text{ and } o(v)\to\theta. \end{aligned} \tag{3.4}$$

At first, we estimate the first part in the limit (3.3).

$$\left\| \frac{\frac{a(u)\bar{x}}{\|a(u)\bar{x}+o(u)\|}-\frac{a(v)\bar{x}}{\|a(v)\bar{x}+o(v)\|}}{\|u-v\|} \right\|$$

$$= \frac{\left\| \frac{a(u)\bar{x}}{\|a(u)\bar{x}+o(u)\|}-\frac{a(v)\bar{x}}{\|a(v)\bar{x}+o(v)\|} \right\|}{\|u-v\|}$$

$$= \frac{\|\bar{x}\|\left| \frac{a(u)}{\|a(u)\bar{x}+o(u)\|}-\frac{a(v)}{\|a(v)\bar{x}+o(v)\|} \right|}{\|u-v\|}$$

$$= \frac{\|\bar{x}\|\left| \frac{a(u)\|a(v)\bar{x}+o(v)\|-a(v)\|a(u)\bar{x}+o(u)\|}{\|a(u)\bar{x}+o(u)\|\|a(v)\bar{x}+o(v)\|} \right|}{\|u-v\|}$$

$$= \frac{\|\bar{x}\|\left| \frac{(a(u)\|a(v)\bar{x}+o(v)\|)^2-(a(v)\|a(u)\bar{x}+o(u)\|)^2}{\|a(u)\bar{x}+o(u)\|\|a(v)\bar{x}+o(v)\|(a(u)\|a(v)\bar{x}+o(v)\|+a(v)\|a(u)\bar{x}+o(u)\|)} \right|}{\|u-v\|}$$

$$= \frac{\|\bar{x}\|\left| \frac{a(u)^2(a(v)^2\|\bar{x}\|^2+\|o(v)\|^2)-a(v)^2(a(u)^2\|\bar{x}\|^2+\|o(u)\|^2)}{\|a(u)\bar{x}+o(u)\|\|a(v)\bar{x}+o(v)\|(a(u)\|a(v)\bar{x}+o(v)\|+a(v)\|a(u)\bar{x}+o(u)\|)} \right|}{\|(a(u)\bar{x}+o(u))-(a(v)\bar{x}+o(v))\|}$$

$$= \frac{\|\bar{x}\| \left| \frac{a(u)^2 \|o(v)\|^2 - a(v)^2 \|o(u)\|^2}{\|a(u)\bar{x}+o(u)\| \, \|a(v)\bar{x}+o(v)\| \, (a(u)\|a(v)\bar{x}+o(v)\| + a(v)\|a(u)\bar{x}+o(u)\|)} \right|}{\|(a(u)-a(v))\bar{x}+(o(u)-o(v))\|}$$

$$= \frac{\|\bar{x}\| \left| \frac{a(u)^2\|o(v)\|^2 - a(u)^2\|o(u)\|^2 + a(u)^2\|o(u)\|^2 - a(v)^2\|o(u)\|^2}{\|a(u)\bar{x}+o(u)\| \, \|a(v)\bar{x}+o(v)\| \, (a(u)\|a(v)\bar{x}+o(v)\| + a(v)\|a(u)\bar{x}+o(u)\|)} \right|}{\sqrt{\big(a(u)-a(v)\big)^2 \|\bar{x}\|^2 + \|o(u)-o(v)\|^2}}$$

$$= \frac{\|\bar{x}\| \left| \frac{a(u)^2(\|o(v)\|+\|o(u)\|)(\|o(v)\|-\|o(u)\|) + \|o(u)\|^2(a(u)+a(v))(a(u)-a(v))}{\|a(u)\bar{x}+o(u)\| \, \|a(v)\bar{x}+o(v)\| \, (a(u)\|a(v)\bar{x}+o(v)\| + a(v)\|a(u)\bar{x}+o(u)\|)} \right|}{\sqrt{\big(a(u)-a(v)\big)^2 \|\bar{x}\|^2 + \|o(u)-o(v)\|^2}}$$

$$\leq \frac{\|\bar{x}\| \left| \frac{a(u)^2(\|o(v)\|+\|o(u)\|)(\|o(v)\|-\|o(u)\|)}{\|a(u)\bar{x}+o(u)\| \, \|a(v)\bar{x}+o(v)\| \, (a(u)\|a(v)\bar{x}+o(v)\| + a(v)\|a(u)\bar{x}+o(u)\|)} \right|}{\sqrt{\big(a(u)-a(v)\big)^2 \|\bar{x}\|^2 + \|o(u)-o(v)\|^2}}$$

$$+ \frac{\|\bar{x}\| \left| \frac{\|o(u)\|^2(a(u)+a(v))(a(u)-a(v))}{\|a(u)\bar{x}+o(u)\| \, \|a(v)\bar{x}+o(v)\| \, (a(u)\|a(v)\bar{x}+o(v)\| + a(v)\|a(u)\bar{x}+o(u)\|)} \right|}{\sqrt{\big(a(u)-a(v)\big)^2 \|\bar{x}\|^2 + \|o(u)-o(v)\|^2}}$$

$$\leq \frac{\|\bar{x}\| \, \frac{a(u)^2(\|o(v)\|+\|o(u)\|) \, | \|o(v)\|-\|o(u)\| |}{\|a(u)\bar{x}+o(u)\| \, \|a(v)\bar{x}+o(v)\| \, (a(u)\|a(v)\bar{x}+o(v)\| + a(v)\|a(u)\bar{x}+o(u)\|)}}{\|o(u)-o(v)\|}$$

$$+ \frac{\|\bar{x}\| \, \frac{\|o(u)\|^2(a(u)+a(v)) \, |a(u)-a(v)|}{\|a(u)\bar{x}+o(u)\| \, \|a(v)\bar{x}+o(v)\| \, (a(u)\|a(v)\bar{x}+o(v)\| + a(v)\|a(u)\bar{x}+o(u)\|)}}{|a(u)-a(v)| \, \|\bar{x}\|}$$

$$\leq \|\bar{x}\| \, \frac{a(u)^2(\|o(v)\|+\|o(u)\|)}{\|a(u)\bar{x}+o(u)\| \, \|a(v)\bar{x}+o(v)\| \, (a(u)\|a(v)\bar{x}+o(v)\| + a(v)\|a(u)\bar{x}+o(u)\|)}$$

$$+ \frac{\|o(u)\|^2(a(u)+a(v))}{\|a(u)\bar{x}+o(u)\| \, \|a(v)\bar{x}+o(v)\| \, (a(u)\|a(v)\bar{x}+o(v)\| + a(v)\|a(u)\bar{x}+o(u)\|)}$$

$$\longrightarrow \frac{0}{2\|\bar{x}\|^2} + \frac{0}{2\|\bar{x}\|^2} = 0, \text{ as } u \to \bar{x} \text{ and } v \to \bar{x}. \tag{3.5}$$

Where, when we take the limit, we applied (3.4). Then, we estimate the second part in (3.3).

$$\frac{\left\| \frac{o(u)}{\|a(u)\bar{x}+o(u)\|} - \frac{o(u)}{\|\bar{x}\|} - \left(\frac{o(v)}{\|a(v)\bar{x}+o(v)\|} - \frac{o(v)}{\|\bar{x}\|}\right) \right\|}{\|u-v\|}$$

$$= \frac{\left\| \frac{o(u)(\|\bar{x}\| - \|a(u)\bar{x}+o(u)\|)}{\|\bar{x}\| \, \|a(u)\bar{x}+o(u)\|} - \left(\frac{o(v)(\|\bar{x}\| - \|a(v)\bar{x}+o(v)\|)}{\|\bar{x}\| \, \|a(v)\bar{x}+o(v)\|}\right) \right\|}{\|u-v\|}$$

$$= \frac{\frac{1}{\|\bar{x}\|} \left\| \frac{o(u)(\|a(u)\bar{x}+o(u)\|^2 - \|\bar{x}\|^2)}{\|a(u)\bar{x}+o(u)\|(\|a(u)\bar{x}+o(u)\|+\|\bar{x}\|)} - \frac{o(v)(\|a(v)\bar{x}+o(v)\|^2 - \|\bar{x}\|^2)}{\|a(v)\bar{x}+o(v)\|(\|a(v)\bar{x}+o(v)\|+\|\bar{x}\|)} \right\|}{\|u-v\|}$$

$$= \frac{\frac{1}{\|\bar{x}\|} \left\| \frac{o(u)((a(u)^2-1)\|\bar{x}\|^2 + \|o(u)\|^2)}{\|a(u)\bar{x}+o(u)\|(\|a(u)\bar{x}+o(u)\|+\|\bar{x}\|)} - \frac{o(v)((a(v)^2-1)\|\bar{x}\|^2 + \|o(v)\|^2)}{\|a(v)\bar{x}+o(v)\|(\|a(v)\bar{x}+o(v)\|+\|\bar{x}\|)} \right\|}{\|u-v\|}$$

$$\leq \frac{\frac{1}{\|\bar{x}\|}\left\|\frac{o(u)((a(u)^2-1)\|\bar{x}\|^2+\|o(u)\|^2)}{\|a(u)\bar{x}+o(u)\|(\|a(u)\bar{x}+o(u)\|+\|\bar{x}\|)} - \frac{o(v)((a(v)^2-1)\|\bar{x}\|^2+\|o(v)\|^2)}{\|a(u)\bar{x}+o(u)\|(\|a(u)\bar{x}+o(u)\|+\|\bar{x}\|)}\right\|}{\sqrt{(a(u)-a(v))^2\|\bar{x}\|^2+\|o(u)-o(v)\|^2}}$$

$$+ \frac{\frac{1}{\|\bar{x}\|}\left\|\frac{o(v)((a(v)^2-1)\|\bar{x}\|^2+\|o(v)\|^2)}{\|a(u)\bar{x}+o(u)\|(\|a(u)\bar{x}+o(u)\|+\|\bar{x}\|)} - \frac{o(v)((a(v)^2-1)\|\bar{x}\|^2+\|o(v)\|^2)}{\|a(v)\bar{x}+o(v)\|(\|a(v)\bar{x}+o(v)\|+\|\bar{x}\|)}\right\|}{\sqrt{(a(u)-a(v))^2\|\bar{x}\|^2+\|o(u)-o(v)\|^2}}. \tag{3.6}$$

We estimate the first part in (3.6).

$$\frac{\frac{1}{\|\bar{x}\|}\left\|\frac{o(u)((a(u)^2-1)\|\bar{x}\|^2+\|o(u)\|^2)}{\|a(u)\bar{x}+o(u)\|(\|a(u)\bar{x}+o(u)\|+\|\bar{x}\|)} - \frac{o(v)((a(v)^2-1)\|\bar{x}\|^2+\|o(v)\|^2)}{\|a(u)\bar{x}+o(u)\|(\|a(u)\bar{x}+o(u)\|+\|\bar{x}\|)}\right\|}{\sqrt{(a(u)-a(v))^2\|\bar{x}\|^2+\|o(u)-o(v)\|^2}}$$

$$= \frac{\frac{1}{\|\bar{x}\|\|a(u)\bar{x}+o(u)\|(\|a(u)\bar{x}+o(u)\|+\|\bar{x}\|)}\left\|o(u)((a(u)^2-1)\|\bar{x}\|^2+\|o(u)\|^2)-o(v)((a(v)^2-1)\|\bar{x}\|^2+\|o(v)\|^2)\right\|}{\sqrt{(a(u)-a(v))^2\|\bar{x}\|^2+\|o(u)-o(v)\|^2}} \tag{3.7}$$

For the simplicity, we write

$$b(u,v) := \|\bar{x}\|\|a(u)\bar{x}+o(u)\|(\|a(u)\bar{x}+o(u)\|+\|\bar{x}\|).$$

Taking limit, we have

$$\lim_{\substack{u\to\bar{x},v\to\bar{x}\\u,v\in\mathbb{B}(\bar{x},p)}} b(u,v)$$

$$= \lim_{\substack{u\to\bar{x},v\to\bar{x}\\u,v\in\mathbb{B}(\bar{x},p)}} \left(\|\bar{x}\|\|a(u)\bar{x}+o(u)\|(\|a(u)\bar{x}+o(u)\|+\|\bar{x}\|)\right)$$

$$= 2\|\bar{x}\|^3.$$

Then, (3.7) can be rewritten as

$$\frac{\frac{1}{\|\bar{x}\|\|a(u)\bar{x}+o(u)\|(\|a(u)\bar{x}+o(u)\|+\|\bar{x}\|)}\left\|o(u)((a(u)^2-1)\|\bar{x}\|^2+\|o(u)\|^2)-o(v)((a(v)^2-1)\|\bar{x}\|^2+\|o(v)\|^2)\right\|}{\sqrt{(a(u)-a(v))^2\|\bar{x}\|^2+\|o(u)-o(v)\|^2}}$$

$$= \frac{\frac{1}{b(u,v)}\left\|o(u)((a(u)^2-1)\|\bar{x}\|^2+\|o(u)\|^2)-o(v)((a(v)^2-1)\|\bar{x}\|^2+\|o(v)\|^2)\right\|}{\sqrt{(a(u)-a(v))^2\|\bar{x}\|^2+\|o(u)-o(v)\|^2}}. \tag{3.7}$$

We calculate the terms inside the norm in the numerator of (3.7).

$$o(u)((a(u)^2-1)\|\bar{x}\|^2+\|o(u)\|^2)-o(v)((a(v)^2-1)\|\bar{x}\|^2+\|o(v)\|^2)$$

$$= (o(u)(a(u)^2-1)-o(v)(a(v)^2-1)))\|\bar{x}\|^2 + o(u)\|o(u)\|^2 - o(v)\|o(v)\|^2$$

$$= (o(u)(a(u)^2-1) - o(u)(a(v)^2-1) + o(u)(a(v)^2-1) - o(v)(a(v)^2-1))\|\bar{x}\|^2$$

$$+ o(u)\|o(u)\|^2 - o(v)\|o(u)\|^2 + o(v)\|o(u)\|^2 - o(v)\|o(v)\|^2$$

$$= (o(u)((a(u)^2 - 1) - (a(v)^2 - 1)) + (o(u) - o(v))(a(v)^2 - 1))\|\bar{x}\|^2$$

$$+ (o(u) - o(v))\|o(u)\|^2 + o(v)(\|o(u)\| - \|o(v)\|)(\|o(u)\| + \|o(v)\|)$$

$$= (o(u)(a(u)^2 - a(v)^2) + (o(u) - o(v))(a(v)^2 - 1))\|\bar{x}\|^2$$

$$+ (o(u) - o(v))\|o(u)\|^2 + o(v)(\|o(u)\| - \|o(v)\|)(\|o(u)\| + \|o(v)\|)$$

$$= (o(u)(a(u) - a(v))(a(u) + a(v)) + (o(u) - o(v))(a(v)^2 - 1))\|\bar{x}\|^2$$

$$+ (o(u) - o(v))\|o(u)\|^2 + o(v)(\|o(u)\| - \|o(v)\|)(\|o(u)\| + \|o(v)\|)$$

Then, (3.7) is estimated by the following 4 parts.

$$\frac{1}{b(u,v)} \frac{\| o(u)((a(u)^2-1)\|\bar{x}\|^2 + \|o(u)\|^2) - o(v)((a(v)^2-1)\|\bar{x}\|^2 + \|o(v)\|^2) \|}{\sqrt{(a(u)-a(v))^2\|\bar{x}\|^2 + \|o(u)-o(v)\|^2}}$$

$$\leq \frac{1}{b(u,v)} \frac{\| o(u)(a(u)-a(v))(a(u)+a(v))\|\bar{x}\|^2 \|}{\sqrt{(a(u)-a(v))^2\|\bar{x}\|^2 + \|o(u)-o(v)\|^2}} + \frac{1}{b(u,v)} \frac{\| (o(u)-o(v))(a(v)^2-1)\|\bar{x}\|^2 \|}{\sqrt{(a(u)-a(v))^2\|\bar{x}\|^2 + \|o(u)-o(v)\|^2}}$$

$$+ \frac{1}{b(u,v)} \frac{\| (o(u)-o(v))\|o(u)\|^2 \|}{\sqrt{(a(u)-a(v))^2\|\bar{x}\|^2 + \|o(u)-o(v)\|^2}} + \frac{1}{b(u,v)} \frac{\| o(v)(\|o(u)\|-\|o(v)\|)(\|o(u)\|+\|o(v)\|)\|}{\sqrt{(a(u)-a(v))^2\|\bar{x}\|^2 + \|o(u)-o(v)\|^2}}$$

$$= \frac{\|o(u)\|\|\bar{x}\|^2}{b(u,v)} \frac{|(a(u)-a(v))(a(u)+a(v))|}{\sqrt{(a(u)-a(v))^2\|\bar{x}\|^2 + \|o(u)-o(v)\|^2}} + \frac{\|\bar{x}\|^2}{b(u,v)} \frac{|a(v)^2-1|\|o(u)-o(v)\|}{\sqrt{(a(u)-a(v))^2\|\bar{x}\|^2 + \|o(u)-o(v)\|^2}}$$

$$+ \frac{\|o(u)\|^2}{b(u,v)} \frac{\|o(u)-o(v)\|}{\sqrt{(a(u)-a(v))^2\|\bar{x}\|^2 + \|o(u)-o(v)\|^2}} + \frac{\|o(v)\|}{b(u,v)} \frac{|\|o(u)\|-\|o(v)\||(\|o(u)\|+\|o(v)\|)}{\sqrt{(a(u)-a(v))^2\|\bar{x}\|^2 + \|o(u)-o(v)\|^2}}$$

$$\leq \frac{\|o(u)\|\|\bar{x}\|^2}{b(u,v)} \frac{|(a(u)-a(v))(a(u)+a(v))|}{|(a(u)-a(v))|\|\bar{x}\|} + \frac{\|\bar{x}\|^2}{b(u,v)} \frac{|a(v)^2-1|\|o(u)-o(v)\|}{\|o(u)-o(v)\|}$$

$$+ \frac{\|o(u)\|^2}{b(u,v)} \frac{\|o(u)-o(v)\|}{\|o(u)-o(v)\|} + \frac{\|o(v)\|}{b(u,v)} \frac{|\|o(u)\|-\|o(v)\||(\|o(u)\|+\|o(v)\|)}{\|o(u)-o(v)\|}$$

$$\leq \frac{\|o(u)\|\|a(u)+a(v)\|\|\bar{x}\|}{b(u,v)} + \frac{\|\bar{x}\|^2|a(v)^2-1|}{b(u,v)} + \frac{\|o(u)\|^2}{b(u,v)} + \frac{\|o(v)\|(\|o(u)\|+\|o(v)\|)}{b(u,v)}$$

$$\to 0, \text{ as } u \to \bar{x} \text{ and } v \to \bar{x}. \tag{3.8}$$

Where, we applied again (3.4) and $b(u,v) \to 2\|\bar{x}\|^3$, as $u \to \bar{x}$ and $v \to \bar{x}$.

Next, we estimate the second term in (3.6).

$$\frac{\frac{1}{\|\bar{x}\|}\left\|\frac{o(v)((a(v)^2-1)\|\bar{x}\|^2+\|o(v)\|^2)}{\|a(u)\bar{x}+o(u)\|(\|a(u)\bar{x}+o(u)\|+\|\bar{x}\|)}-\frac{o(v)((a(v)^2-1)\|\bar{x}\|^2+\|o(v)\|^2)}{\|a(v)\bar{x}+o(v)\|(\|a(v)\bar{x}+o(v)\|+\|\bar{x}\|)}\right\|}{\sqrt{(a(u)-a(v))^2\|\bar{x}\|^2+\|o(u)-o(v)\|^2}}$$

$$=\frac{\frac{\|o(v)\|(|(a(v)^2-1|\|\bar{x}\|^2+\|o(v)\|^2)}{\|\bar{x}\|}\left|\frac{1}{\|a(u)\bar{x}+o(u)\|(\|a(u)\bar{x}+o(u)\|+\|\bar{x}\|)}-\frac{1}{\|a(v)\bar{x}+o(v)\|(\|a(v)\bar{x}+o(v)\|+\|\bar{x}\|)}\right|}{\sqrt{(a(u)-a(v))^2\|\bar{x}\|^2+\|o(u)-o(v)\|^2}}$$

$$=\frac{\frac{\|o(v)\|(|(a(v)^2-1|\|\bar{x}\|^2+\|o(v)\|^2)}{\|\bar{x}\|}\left|\frac{\|a(u)\bar{x}+o(u)\|(\|a(u)\bar{x}+o(u)\|+\|\bar{x}\|)-\|a(v)\bar{x}+o(v)\|(\|a(v)\bar{x}+o(v)\|+\|\bar{x}\|)}{\|a(u)\bar{x}+o(u)\|(\|a(u)\bar{x}+o(u)\|+\|\bar{x}\|)\|a(v)\bar{x}+o(v)\|(\|a(v)\bar{x}+o(v)\|+\|\bar{x}\|)}\right|}{\sqrt{(a(u)-a(v))^2\|\bar{x}\|^2+\|o(u)-o(v)\|^2}}$$

$$=\frac{\frac{\|o(v)\|\|o(v)\|(|(a(v)^2-1|\|\bar{x}\|^2+\|o(v)\|^2)}{\|a(u)\bar{x}+o(u)\|(\|a(u)\bar{x}+o(u)\|+\|\bar{x}\|)\|a(v)\bar{x}+o(v)\|(\|a(v)\bar{x}+o(v)\|+\|\bar{x}\|)\|\bar{x}\|}\left|\frac{\|a(u)\bar{x}+o(u)\|(\|a(u)\bar{x}+o(u)\|+\|\bar{x}\|)-\|a(v)\bar{x}+o(v)\|(\|a(v)\bar{x}+o(v)\|+\|\bar{x}\|)}{1}\right|}{\sqrt{(a(u)-a(v))^2\|\bar{x}\|^2+\|o(u)-o(v)\|^2}}$$

$$=\frac{\frac{\|o(v)\|\left(\left|(a(v)^2-1\right|\|\bar{x}\|^2+\|o(v)\|^2\right)}{\|a(u)\bar{x}+o(u)\|(\|a(u)\bar{x}+o(u)\|+\|\bar{x}\|)\|a(v)\bar{x}+o(v)\|(\|a(v)\bar{x}+o(v)\|+\|\bar{x}\|)\|\bar{x}\|}\left|\frac{\|a(u)\bar{x}+o(u)\|^2(\|a(u)\bar{x}+o(u)\|+\|\bar{x}\|)^2-\|a(v)\bar{x}+o(v)\|^2(\|a(v)\bar{x}+o(v)\|+\|\bar{x}\|)^2}{\|a(u)\bar{x}+o(u)\|(\|a(u)\bar{x}+o(u)\|+\|\bar{x}\|)+\|a(v)\bar{x}+o(v)\|(\|a(v)\bar{x}+o(v)\|+\|\bar{x}\|)}\right|}{\sqrt{(a(u)-a(v))^2\|\bar{x}\|^2+\|o(u)-o(v)\|^2}} \quad (3.9)$$

For simplicity, we write

$$c(u,v) := \|o(v)\|(|(a(v)^2-1|\|\bar{x}\|^2+\|o(v)\|^2).$$

Taking limit gets

$$\lim_{\substack{u\to\bar{x},v\to\bar{x}\\u,v\in\mathbb{B}(\bar{x},p)}} c(u,v) = \lim_{\substack{u\to\bar{x},v\to\bar{x}\\u,v\in\mathbb{B}(\bar{x},p)}} \|o(v)\|(|(a(v)^2-1|\|\bar{x}\|^2+\|o(v)\|^2) = 0. \quad (3.10)$$

And, we writ

$$d(u,v) := [\|a(u)\bar{x}+o(u)\|(\|a(u)\bar{x}+o(u)\|+\|\bar{x}\|)\|a(v)\bar{x}+o(v)\|(\|a(v)\bar{x}+o(v)\|+\|\bar{x}\|)\|\bar{x}\|]$$
$$[\|a(u)\bar{x}+o(u)\|(\|a(u)\bar{x}+o(u)\|+\|\bar{x}\|)+\|a(v)\bar{x}+o(v)\|(\|a(v)\bar{x}+o(v)\|+\|\bar{x}\|)].$$

It satisfies

$$\lim_{\substack{u\to\bar{x},v\to\bar{x}\\u,v\in\mathbb{B}(\bar{x},p)}} d(u,v) = 4\|\bar{x}\|^5 \times 4\|\bar{x}\|^2 = 16\|\bar{x}\|^7. \quad (3.11)$$

Then, the second term of (3.6), which is (3.9), is rewritten as

$$\frac{c(u,v)}{d(u,v)}\frac{\left|\|a(u)\bar{x}+o(u)\|^2(\|a(u)\bar{x}+o(u)\|+\|\bar{x}\|)^2-\|a(v)\bar{x}+o(v)\|^2(\|a(v)\bar{x}+o(v)\|+\|\bar{x}\|)^2\right|}{\sqrt{(a(u)-a(v))^2\|\bar{x}\|^2+\|o(u)-o(v)\|^2}}. \quad (3.9)$$

We estimate the crucial factor in (3.9), which is the second factor.

$$\frac{\left|\|a(u)\bar{x}+o(u)\|^2(\|a(u)\bar{x}+o(u)\|+\|\bar{x}\|)^2-\|a(v)\bar{x}+o(v)\|^2(\|a(v)\bar{x}+o(v)\|+\|\bar{x}\|)^2\right|}{\sqrt{(a(u)-a(v))^2\|\bar{x}\|^2+\|o(u)-o(v)\|^2}}$$

$$= \frac{\left|(a(u)^2\|\bar{x}\|^2+\|o(u)\|^2)(\|a(u)\bar{x}+o(u)\|^2+2\|a(u)\bar{x}+o(u)\|\|\bar{x}\|+\|\bar{x}\|^2)\right.}{\sqrt{(a(u)-a(v))^2\|\bar{x}\|^2+\|o(u)-o(v)\|^2}}$$

$$-\frac{(a(v)^2\|\bar{x}\|^2+\|o(v)\|^2)(\|a(v)\bar{x}+o(v)\|^2+2\|a(v)\bar{x}+o(v)\|\|\bar{x}\|+\|\bar{x}\|^2)\left.\right|}{\sqrt{(a(u)-a(v))^2\|\bar{x}\|^2+\|o(u)-o(v)\|^2}}$$

$$= \frac{\left|(a(u)^2\|\bar{x}\|^2+\|o(u)\|^2)(a(u)^2\|\bar{x}\|^2+\|o(u)\|^2+2\|a(u)\bar{x}+o(u)\|\|\bar{x}\|+\|\bar{x}\|^2)\right.}{\sqrt{(a(u)-a(v))^2\|\bar{x}\|^2+\|o(u)-o(v)\|^2}}$$

$$-\frac{(a(v)^2\|\bar{x}\|^2+\|o(v)\|^2)(a(u)^2\|\bar{x}\|^2+\|o(u)\|^2+2\|a(u)\bar{x}+o(u)\|\|\bar{x}\|+\|\bar{x}\|^2)}{\sqrt{(a(u)-a(v))^2\|\bar{x}\|^2+\|o(u)-o(v)\|^2}}$$

$$+\frac{(a(v)^2\|\bar{x}\|^2+\|o(v)\|^2)(a(u)^2\|\bar{x}\|^2+\|o(u)\|^2+2\|a(u)\bar{x}+o(u)\|\|\bar{x}\|+\|\bar{x}\|^2)}{\sqrt{(a(u)-a(v))^2\|\bar{x}\|^2+\|o(u)-o(v)\|^2}}$$

$$-\frac{(a(v)^2\|\bar{x}\|^2+\|o(v)\|^2)(a(v)^2\|\bar{x}\|^2+\|o(v)\|^2+2\|a(v)\bar{x}+o(v)\|\|\bar{x}\|+\|\bar{x}\|^2)\left.\right|}{\sqrt{(a(u)-a(v))^2\|\bar{x}\|^2+\|o(u)-o(v)\|^2}}$$

$$= \frac{\left|\left((a(u)^2\|\bar{x}\|^2+\|o(u)\|^2)-(a(v)^2\|\bar{x}\|^2+\|o(v)\|^2)\right)(a(u)^2\|\bar{x}\|^2+\|o(u)\|^2+2\|a(u)\bar{x}+o(u)\|\|\bar{x}\|+\|\bar{x}\|^2)\right.}{\sqrt{(a(u)-a(v))^2\|\bar{x}\|^2+\|o(u)-o(v)\|^2}}$$

$$+\frac{(a(v)^2\|\bar{x}\|^2+\|o(v)\|^2)\left((a(u)^2\|\bar{x}\|^2+\|o(u)\|^2+2\|a(u)\bar{x}+o(u)\|\|\bar{x}\|+\|\bar{x}\|^2)-(a(v)^2\|\bar{x}\|^2+\|o(v)\|^2+2\|a(v)\bar{x}+o(v)\|\|\bar{x}\|+\|\bar{x}\|^2)\right)}{\sqrt{(a(u)-a(v))^2\|\bar{x}\|^2+\|o(u)-o(v)\|^2}}$$

$$= \frac{\left|\left((a(u)^2-a(v)^2)\|\bar{x}\|^2+\|o(u)\|^2-\|o(v)\|^2\right)(a(u)^2\|\bar{x}\|^2+\|o(u)\|^2+2\|a(u)\bar{x}+o(u)\|\|\bar{x}\|+\|\bar{x}\|^2)\right.}{\sqrt{(a(u)-a(v))^2\|\bar{x}\|^2+\|o(u)-o(v)\|^2}}$$

$$+\frac{(a(v)^2\|\bar{x}\|^2+\|o(v)\|^2)\left((a(u)^2\|\bar{x}\|^2+\|o(u)\|^2+2\|a(u)\bar{x}+o(u)\|\|\bar{x}\|)-(a(v)^2\|\bar{x}\|^2+\|o(v)\|^2+2\|a(v)\bar{x}+o(v)\|\|\bar{x}\|)\right)\left.\right|}{\sqrt{(a(u)-a(v))^2\|\bar{x}\|^2+\|o(u)-o(v)\|^2}}$$

$$= \frac{\left|\left((a(u)-a(v))(a(u)+a(v))\|\bar{x}\|^2+(\|o(u)\|-\|o(v)\|)(\|o(u)\|+\|o(v)\|)\right)(a(u)^2\|\bar{x}\|^2+\|o(u)\|^2+2\|a(u)\bar{x}+o(u)\|\|\bar{x}\|+\|\bar{x}\|^2)\right.}{\sqrt{(a(u)-a(v))^2\|\bar{x}\|^2+\|o(u)-o(v)\|^2}}$$

$$+\frac{(a(v)^2\|\bar{x}\|^2+\|o(v)\|^2)\left(((a(u)-a(v))(a(u)+a(v))\|\bar{x}\|^2+(\|o(u)\|-\|o(v)\|)(\|o(u)\|+\|o(v)\|)+2(\|a(u)\bar{x}+o(u)\|-\|a(v)\bar{x}+o(v)\|)\|\bar{x}\|)\right)\left.\right|}{\sqrt{(a(u)-a(v))^2\|\bar{x}\|^2+\|o(u)-o(v)\|^2}}$$

$$\leq \frac{\left|(a(u)-a(v))(a(u)+a(v))\|\bar{x}\|^2(a(u)^2\|\bar{x}\|^2+\|o(u)\|^2+2\|a(u)\bar{x}+o(u)\|\|\bar{x}\|+\|\bar{x}\|^2)\right|}{\sqrt{(a(u)-a(v))^2\|\bar{x}\|^2+\|o(u)-o(v)\|^2}}$$

$$+\frac{\left|(\|o(u)\|-\|o(v)\|)(\|o(u)\|+\|o(v)\|)(a(u)^2\|\bar{x}\|^2+\|o(u)\|^2+2\|a(u)\bar{x}+o(u)\|\|\bar{x}\|+\|\bar{x}\|^2)\right|}{\sqrt{(a(u)-a(v))^2\|\bar{x}\|^2+\|o(u)-o(v)\|^2}}$$

$$+ \frac{\left|\left(a(v)^2\|\bar{x}\|^2+\|o(v)\|^2\right)\left(a(u)-a(v)\right)\left(a(u)+a(v)\right)\|\bar{x}\|^2\right|}{\sqrt{\left(a(u)-a(v)\right)^2\|\bar{x}\|^2+\|o(u)-o(v)\|^2}}$$

$$+ \frac{\left|\left(a(v)^2\|\bar{x}\|^2+\|o(v)\|^2\right)(\|o(u)\|-\|o(v)\|)(\|o(u)\|+\|o(v)\|)\right|}{\sqrt{\left(a(u)-a(v)\right)^2\|\bar{x}\|^2+\|o(u)-o(v)\|^2}}$$

$$+ \frac{\left|2\left(a(v)^2\|\bar{x}\|^2+\|o(v)\|^2\right)\right\|\left(a(u)\bar{x}+o(u)\right)-\left(a(v)\bar{x}+o(v)\right)\|\|\bar{x}\|}{\sqrt{\left(a(u)-a(v)\right)^2\|\bar{x}\|^2+\|o(u)-o(v)\|^2}}$$

$$\leq \frac{|a(u)-a(v)||a(u)+a(v)|\|\bar{x}\|^2\left(a(u)^2\|\bar{x}\|^2+\|o(u)\|^2+2\|a(u)\bar{x}+o(u)\|\|\bar{x}\|+\|\bar{x}\|^2\right)}{|a(u)-a(v)|\|\bar{x}\|}$$

$$+ \frac{\|\|o(u)\|-\|o(v)\|\|(\|o(u)\|+\|o(v)\|)\left(a(u)^2\|\bar{x}\|^2+\|o(u)\|^2+2\|a(u)\bar{x}+o(u)\|\|\bar{x}\|+\|\bar{x}\|^2\right)}{\|o(u)-o(v)\|}$$

$$+ \frac{\left(a(v)^2\|\bar{x}\|^2+\|o(v)\|^2\right)|a(u)-a(v)||a(u)+a(v)|\|\bar{x}\|^2}{|a(u)-a(v)|\|\bar{x}\|}$$

$$+ \frac{\left(a(v)^2\|\bar{x}\|^2+\|o(v)\|^2\right)\|\|o(u)\|-\|o(v)\|\|(\|o(u)\|+\|o(v)\|)}{\|o(u)-o(v)\|}$$

$$+ \frac{2\left(a(v)^2\|\bar{x}\|^2+\|o(v)\|^2\right)\|\left(a(u)\bar{x}+o(u)\right)-\left(a(v)\bar{x}+o(v)\right)\|\|\bar{x}\|}{\sqrt{\left(a(u)-a(v)\right)^2\|\bar{x}\|^2+\|o(u)-o(v)\|^2}}$$

$$\leq \frac{|a(u)+a(v)|\|\bar{x}\|^2\left(a(u)^2\|\bar{x}\|^2+\|o(u)\|^2+2\|a(u)\bar{x}+o(u)\|\|\bar{x}\|+\|\bar{x}\|^2\right)}{\|\bar{x}\|}$$

$$+ \frac{\|o(u)-o(v)\|(\|o(u)\|+\|o(v)\|)\left(a(u)^2\|\bar{x}\|^2+\|o(u)\|^2+2\|a(u)\bar{x}+o(u)\|\|\bar{x}\|+\|\bar{x}\|^2\right)}{\|o(u)-o(v)\|}$$

$$+ \left(a(v)^2\|\bar{x}\|^2+\|o(v)\|^2\right)|a(u)+a(v)|\|\bar{x}\|$$

$$+ \frac{\left(a(v)^2\|\bar{x}\|^2+\|o(v)\|^2\right)\|o(u)-o(v)\|(\|o(u)\|+\|o(v)\|)}{\|o(u)-o(v)\|}$$

$$+ \frac{2\left(a(v)^2\|\bar{x}\|^2+\|o(v)\|^2\right)(|a(u)-a(v)|\|\bar{x}\|+\|o(u)-o(v)\|)\|\bar{x}\|}{\sqrt{\left(a(u)-a(v)\right)^2\|\bar{x}\|^2+\|o(u)-o(v)\|^2}}$$

$$\leq \frac{|a(u)+a(v)|\|\bar{x}\|^2\left(a(u)^2\|\bar{x}\|^2+\|o(u)\|^2+2\|a(u)\bar{x}+o(u)\|\|\bar{x}\|+\|\bar{x}\|^2\right)}{\|\bar{x}\|}$$

$$+(\|o(u)\|+\|o(v)\|)\left(a(u)^2\|\bar{x}\|^2+\|o(u)\|^2+2\|a(u)\bar{x}+o(u)\|\|\bar{x}\|+\|\bar{x}\|^2\right)$$

$$+ \left(a(v)^2\|\bar{x}\|^2+\|o(v)\|^2\right)|a(u)+a(v)|\|\bar{x}\|$$

$$+ \left(a(v)^2\|\bar{x}\|^2+\|o(v)\|^2\right)(\|o(u)\|+\|o(v)\|)$$

$$+ \frac{2\left(a(v)^2\|\bar{x}\|^2+\|o(v)\|^2\right)|a(u)-a(v)|\|\bar{x}\|^2}{\sqrt{\left(a(u)-a(v)\right)^2\|\bar{x}\|^2+\|o(u)-o(v)\|^2}}$$

$$+ \frac{2(a(v)^2\|\bar{x}\|^2+\|o(v)\|^2)\|o(u)-o(v)\|\,\|\bar{x}\|}{\sqrt{(a(u)-a(v))^2\|\bar{x}\|^2+\|o(u)-o(v)\|^2}}$$

$$\leq \frac{|a(u)+a(v)|\|\bar{x}\|^2(a(u)^2\|\bar{x}\|^2+\|o(u)\|^2+2\|a(u)\bar{x}+o(u)\|\|\bar{x}\|+\|\bar{x}\|^2)}{\|\bar{x}\|}$$

$$+(\|o(u)\|+\|o(v)\|)(a(u)^2\|\bar{x}\|^2+\|o(u)\|^2+2\|a(u)\bar{x}+o(u)\|\|\bar{x}\|+\|\bar{x}\|^2)$$

$$+(a(v)^2\|\bar{x}\|^2+\|o(v)\|^2)|a(u)+a(v)|\|\bar{x}\|$$

$$+(a(v)^2\|\bar{x}\|^2+\|o(v)\|^2)(\|o(u)\|+\|o(v)\|)$$

$$+2(a(v)^2\|\bar{x}\|^2+\|o(v)\|^2)\|\bar{x}\|$$

$$+2(a(v)^2\|\bar{x}\|^2+\|o(v)\|^2)\,\|\bar{x}\|$$

$$\longrightarrow 8\|\bar{x}\|^3 + 0 + 2\|\bar{x}\|^3 + 0 + 2\|\bar{x}\|^3 + 2\|\bar{x}\|^3$$

$$= 14\|\bar{x}\|^3, \text{ as } u \to \bar{x} \text{ and } v \to \bar{x}. \tag{3.12}$$

By (3.10), (3.11) and (3.12), taking limit for (3.9), we get the limit of the second term in (3.6).

$$\lim_{\substack{u\to\bar{x},v\to\bar{x}\\u,v\in\mathbb{B}(\bar{x},p)}} \frac{\frac{1}{\|\bar{x}\|}\left\|\frac{o(v)((a(v)^2-1)\|\bar{x}\|^2+\|o(v)\|^2)}{\|a(u)\bar{x}+o(u)\|(\|a(u)\bar{x}+o(u)\|+\|\bar{x}\|)} - \frac{o(v)((a(v)^2-1)\|\bar{x}\|^2+\|o(v)\|^2)}{\|a(v)\bar{x}+o(v)\|(\|a(v)\bar{x}+o(v)\|+\|\bar{x}\|)}\right\|}{\sqrt{(a(u)-a(v))^2\|\bar{x}\|^2+\|o(u)-o(v)\|^2}}$$

$$= \lim_{\substack{u\to\bar{x},v\to\bar{x}\\u,v\in\mathbb{B}(\bar{x},p)}} \frac{c(u,v)}{d(u,v)} \frac{\big|\|a(u)\bar{x}+o(u)\|^2(\|a(u)\bar{x}+o(u)\|+\|\bar{x}\|)^2 - \|a(v)\bar{x}+o(v)\|^2(\|a(v)\bar{x}+o(v)\|+\|\bar{x}\|)^2\big|}{\sqrt{(a(u)-a(v))^2\|\bar{x}\|^2+\|o(u)-o(v)\|^2}}$$

$$= \frac{0}{16\|\bar{x}\|^7} 14\|\bar{x}\|^3$$

$$= 0. \tag{3.13}$$

By (3.5), (3.8) and (3.13), we obtain the limit of (3.3) as below.

$$\lim_{u\to\bar{x},v\to\bar{x}} \frac{P_{r\mathbb{B}}(u)-P_{r\mathbb{B}}(v)-\frac{r}{\|\bar{x}\|}o(u-v)}{\|u-v\|} = \lim_{u\to\bar{x},v\to\bar{x}} \frac{P_{r\mathbb{B}}(a(u)\bar{x}+o(u))-P_{r\mathbb{B}}(a(v)\bar{x}+o(v))-\frac{r}{\|\bar{x}\|}o(u-v)}{\|u-v\|} = \theta.$$

This proves part (ii) of this theorem.

Proof of (iii). Part (I) of (iii) follows from part (iii) in Theorem 4.2 in [10] or Theorem 5.2 in [8]. So, we only prove part (II) of (iii). For an arbitrary given $\bar{x} \in r\mathbb{S}$, assume, by the way of contradiction, that $P_{r\mathbb{B}}$ is Fréchet differentiable at $\bar{x}$. Then, there is a linear continuous mapping $A(\bar{x}): H \to H$, such that

$$\lim_{x\to\bar{x}} \frac{P_{r\mathbb{B}}(x)-P_{r\mathbb{B}}(\bar{x})-A(\bar{x})(x-\bar{x})}{\|x-\bar{x}\|} = \theta.$$

In particular, in the above limit, we take a directional line segment $(1 + \delta)\bar{x}$, for $\delta \downarrow 0$. Since $A(\bar{x})$ is assumed to be linear and continuous, we have

$$\theta = \lim_{\delta \downarrow 0} \frac{P_{r\mathbb{B}}((1+\delta)\bar{x}) - P_{r\mathbb{B}}(\bar{x}) - A(\bar{x})((1+\delta)\bar{x} - \bar{x})}{\|(1+\delta)\bar{x} - \bar{x}\|}$$

$$= \lim_{\delta \downarrow 0} \frac{\frac{r}{\|(1+\delta)\bar{x}\|}(1+\delta)\bar{x} - \bar{x} - \delta A(\bar{x})(\bar{x})}{\delta \|\bar{x}\|}$$

$$= \lim_{\delta \downarrow 0} \frac{\frac{r}{(1+\delta)\|\bar{x}\|}(1+\delta)\bar{x} - \bar{x} - \delta A(\bar{x})(\bar{x})}{\delta \|\bar{x}\|}$$

$$= \lim_{\delta \downarrow 0} \frac{\bar{x} - \bar{x}}{\delta \|\bar{x}\|} - \frac{A(\bar{x})(\bar{x})}{\|\bar{x}\|}$$

$$= \lim_{\delta \downarrow 0} \frac{\theta}{\delta r} - \frac{A(\bar{x})(\bar{x})}{r}$$

$$= -\frac{A(\bar{x})(\bar{x})}{r}.$$

This implies
$$A(\bar{x})(\bar{x}) = \theta. \tag{3.14}$$

Next, we take an opposite directional line segment $x = (1 - \delta)\bar{x}$, for $\delta \downarrow 0$ with $0 < \delta < 1$. Since $\|\bar{x}\| = r$, it follows that $(1 - \delta)\bar{x} \in r\mathbb{B}$, for any $\delta$ with $0 < \delta < 1$. By (2.2) and by the assumed linearity of $A(\bar{x})$, we have

$$\theta = \lim_{x \to \bar{x}} \frac{P_{r\mathbb{B}}(x) - P_{r\mathbb{B}}(\bar{x}) - A(\bar{x})(x - \bar{x})}{\|x - \bar{x}\|}$$

$$= \lim_{\delta \downarrow 0, \delta < 1} \frac{(1-\delta)\bar{x} - \bar{x} - A(\bar{x})((1-\delta)\bar{x} - \bar{x})}{\|(1-\delta)\bar{x} - \bar{x}\|}$$

$$= \lim_{\delta \downarrow 0, \delta < 1} \frac{-\delta \bar{x} + \delta A(\bar{x})(\bar{x})}{\delta \|\bar{x}\|}$$

$$= \frac{-\bar{x} + A(\bar{x})(\bar{x})}{r}.$$

This implies
$$A(\bar{x})(\bar{x}) = \bar{x}.$$

This contradicts to (3.14), which proves that $P_{r\mathbb{B}}$ is not Fréchet differentiable at any point $\bar{x} \in r\mathbb{S}$. This theorem is proved. □

Next, we consider some Hilbert spaces with orthonormal bases. Let $\mathbb{N}$ denote the set of all positive integers. Suppose that the considered Hilbert space $H$ has an orthonormal basis $\{e_n : n \in N\}$, in which $N$ is a nonempty subset of $\mathbb{N}$ such that, for any $m, n \in N$, one has

(i) $\langle e_m, e_n \rangle = \begin{cases} 1, & \text{if } m = n, \\ 0, & \text{if } m \neq n. \end{cases}$

(ii) $x = \sum_{n \in N} \langle x, e_n \rangle e_n$, and $\|x\|^2 = \sum_{n \in N} \langle x, e_n \rangle^2$, for every $x \in H$;

(iii) $\langle x, y \rangle = \sum_{n \in N} \langle x, e_n \rangle \langle y, e_n \rangle$, for any $x, y \in H$.

Here, the purpose for us to take a nonempty subset $N$ of $\mathbb{N}$ for the considered Hilbert space is to include all standard Euclidean spaces in the results of Theorem 3.3 and Corollary 3.4.

Let $M$ be an arbitrary nonempty subset of $N$. Let $S(M)$ denote the subspace of $H$ generated by the set $\{e_m : m \in M\}$. Then $H$ has the following orthogonal decomposition

$$H = S(M) \oplus S(N \setminus M).$$

Now we have the following corollary of Theorem 3.3.

**Corollary 3.4.** *Let $H$ be a Hilbert space with an orthonormal basis $\{e_n : n \in N\}$. Let $r > 0$. Then, $P_{r\mathbb{B}}$ has the following properties.*

(i) $P_{r\mathbb{B}}$ *is strictly Fréchet differentiable on $r\mathbb{B}^\circ$. For any $\bar{x} \in r\mathbb{B}^\circ$,*

$$\nabla P_{r\mathbb{B}}(\bar{x})(x) = x, \text{ for every } x \in H.$$

(ii) $P_{r\mathbb{B}}$ *is strictly Fréchet differentiable at every $\bar{x} \in H \setminus r\mathbb{B}$ and for any $x \in H$,*

$$\nabla P_{r\mathbb{B}}(\bar{x})(x) = \frac{r}{\sqrt{\sum_{n \in N} \langle \bar{x}, e_n \rangle^2}} \left( \sum_{n \in N} \langle x, e_n \rangle e_n - \frac{\sum_{n \in N} \langle x, e_n \rangle \langle \bar{x}, e_n \rangle}{\sum_{n \in N} \langle \bar{x}, e_n \rangle^2} \sum_{n \in N} \langle \bar{x}, e_n \rangle e_n \right).$$

(iii) $P_{r\mathbb{B}}$ *is not Fréchet differentiable at any $\bar{x} \in r\mathbb{S}$. That is,*

$$\nabla P_{r\mathbb{B}}(\bar{x}) \text{ does not exist, for any } \bar{x} \in r\mathbb{S}.$$

We give some examples below to demonstrate the results of Theorem 3.3 (or Corollary 3.4), in which the considered Hilbert spaces are just one-, two- or three-dimensional Euclidean spaces.

**Example 3.5.** Let $H = \mathbb{R}$. For any $r > 0$, $r\mathbb{B} = [-r, r]$ and $r\mathbb{S} = \{-r, r\}$. For any real number $\bar{x} \neq 0$, we have

$$\nabla P_{r\mathbb{B}}(\bar{x})(s) = \begin{cases} s, & \text{if } \bar{x} \in (-r, r), \\ 0, & \text{if } \bar{x} \notin [-r, r], \text{ for any } s \in \mathbb{R}. \\ \text{does not exist}, & \text{if } |\bar{x}| = r, \end{cases}$$

**Example 3.6.** Let $H = \mathbb{R}^2 = \{(s, t) : s, t \in \mathbb{R}\}$. For any $r > 0$, $r\mathbb{B}$ is the closed disk in $\mathbb{R}^2$ with radius $r$ and center $\theta := (0, 0)$. Take any point $\bar{x} \in \mathbb{R}^2$ with $\bar{x} \neq \theta$. We consider the following two cases:

(I) $\bar{x} = (a, 0)$ with $a \neq 0$. Then, for any $(s, t) \in \mathbb{R}^2$, we have

$$\nabla P_{r\mathbb{B}}(\bar{x})((s,t)) = \begin{cases} (s,t), & \text{if } a \in (-r,r), \\ \frac{r}{|a|}(0,t), & \text{if } a \notin [-r,r], \\ \text{does not exist}, & \text{if } |a| = r. \end{cases}$$

(II) $\bar{x} = (a, b)$ with $a \neq 0$ and $b \neq 0$. Then,

(i) If $a^2 + b^2 < r^2$, one has

$$\nabla P_{r\mathbb{B}}(\bar{x})((s,t)) = (s,t), \text{ for any } (s, t) \in \mathbb{R}^2;$$

(ii) If $a^2 + b^2 > r^2$, we have

$$\nabla P_{r\mathbb{B}}(\bar{x})((s,t)) = \frac{r}{\sqrt{a^2+b^2}}\left((s,t) - \frac{as+bt}{a^2+b^2}(a,b)\right), \text{ for any } (s, t) \in \mathbb{R}^2;$$

(iii) If $a^2 + b^2 = r^2$, one has

$$\nabla P_{r\mathbb{B}}(\bar{x}) \text{ does not exist};$$

**Example 3.7.** Let $H = \mathbb{R}^3 = \{(s, t, w): s, t, w \in \mathbb{R}\}$. For any $r > 0$, $r\mathbb{B}$ is the closed ball in $\mathbb{R}^3$ with radius $r$ and center $\theta := (0, 0, 0)$. Take any point $\bar{x} \in \mathbb{R}^3$ with $\bar{x} \neq \theta$. We consider the following three cases:

(I) $\bar{x} = (a, 0, 0)$ with $a \neq 0$. we have

$$\nabla P_{r\mathbb{B}}(\bar{x})((s,t,w)) = \begin{cases} (s,t,w), & \text{if } a \in (-r,r), \\ \frac{r}{|a|}(0,t,w), & \text{if } a \notin [-r,r], \\ \text{does not exist}, & \text{if } |a| = r, \end{cases} \text{ for any } (s, t, w) \in \mathbb{R}^3;$$

(II) $\bar{x} = (a, b, 0)$ with $a \neq 0$ and $b \neq 0$. Then,

(i) If $a^2 + b^2 < r^2$, one has

$$\nabla P_{r\mathbb{B}}(\bar{x})((s,t,w)) = (s,t,w), \text{ for any } (s, t, w) \in \mathbb{R}^3;$$

(ii) If $a^2 + b^2 > r^2$, we have

$$\nabla P_{r\mathbb{B}}(\bar{x})((s,t,w)) = \frac{r}{\sqrt{a^2+b^2}}\left((s,t,w) - \frac{as+bt}{a^2+b^2}(a,b,0)\right), \text{ for any } (s, t, w) \in \mathbb{R}^3,$$

(iii) If $a^2 + b^2 = r^2$, one has

$$\nabla P_{r\mathbb{B}}(\bar{x}) \text{ does not exist};$$

(III) $\bar{x} = (a, b, c)$ with $a \neq 0$, $b \neq 0$ and $c \neq 0$. Then

(i) If $a^2 + b^2 + c^2 < r^2$, one has

$$\nabla P_{r\mathbb{B}}(\bar{x})((s,t,w)) = (s,t,w), \text{ for any } (s,t,w) \in \mathbb{R}^3;$$

(ii) If $a^2 + b^2 + c^2 > r^2$, we have

$$\nabla P_{r\mathbb{B}}(\bar{x})((s,t)) = \frac{r}{\sqrt{a^2+b^2+c^2}}\left((s,t,w) - \frac{as+bt+cw}{a^2+b^2+c^2}(a,b,c)\right), \text{ for any } (s,t,w) \in \mathbb{R}^3;$$

(iii) If $a^2 + b^2 + c^2 = r^2$, one has

$$\nabla P_{r\mathbb{B}}(\bar{x}) \text{ does not exist.}$$

## 4. Strict Fréchet differentiability of the metric projection onto the positive cone in $\mathbb{R}^n$

As usual, let $\mathbb{R}^n$ denote the $n$-dimensional Euclidean space with the origin $\theta = (0, 0, \ldots, 0)$. We define

$$\Delta\mathbb{R}^n = \{x = (x_1, x_2, \ldots, x_n) \in \mathbb{R}^n : x_k = 0, \text{ for at least one } k \in \{1, 2, \ldots, n\}\}.$$

Let $K$ denote the positive cone of $\mathbb{R}^n$, which is defined by

$$K = \{x = (x_1, x_2, \ldots, x_n) \in \mathbb{R}^n : x_i \geq 0, i = 1, 2, \ldots, n\}.$$

$K$ is a pointed closed and convex cone in $\mathbb{R}^n$. The interior of $K$ is denoted by $K^o$, which is a nonempty subset of $K$ satisfying

$$K^o = \{x = (x_1, x_2, \ldots, x_n) \in \mathbb{R}^n : x_i > 0, i = 1, 2, \ldots, n\}.$$

The boundary of $K$ is denoted by $\partial K$ such that

$$\partial K = \{x = (x_1, x_2, \ldots, x_n) \in K : x_j = 0, \text{ for at least one } j \in \{1, 2, \ldots, n\}\}.$$

Define

$$\widehat{K} = \{x = (x_1, x_2, \ldots, x_n) \in \mathbb{R}^n : |x_i| > 0, \text{ for all } i \in \{1, 2, \ldots, n\}$$
$$\text{and there are at least one pair } j, k \in \{1, 2, \ldots, n\} \text{ with } x_j x_k < 0\}.$$

Then, we see that $\partial K \subseteq \Delta\mathbb{R}^n$. The negative cone of $K$ is $-K$ satisfying $K \cap (-K) = \{\theta\}$. $-K$ is also a pointed closed and convex cone in $\mathbb{R}^n$. We can similarly define the interior $(-K)^o$. The boundary of $-K$ is denoted by $\partial(-K)$. For any $x = (x_1, x_2, \ldots, x_n) \in \mathbb{R}^n$, we define three subsets of the set $\{1, 2, \ldots, n\}$ with respect to the given $x$ by

$$x^+ = \{i \in \{1, 2, \ldots, n\} : x_i > 0\},$$

$$x^- = \{i \in \{1, 2, \ldots, n\} : x_i < 0\},$$

and
$$\dot{x} = \{i \in \{1, 2, \ldots, n\} : x_i = 0\}.$$

Then, for any $x = (x_1, x_2, \ldots, x_n) \in \mathbb{R}^n$, we have

(a) $x^+ \cup x^- \cup \dot{x} = \{1, 2, \ldots, n\}$;
(b) $x \in K \iff x^- = \emptyset$;
(c) $x \in K^o \iff x^- = \emptyset$ and $\dot{x} = \emptyset$;
(d) $x \in \partial K \iff x^- = \emptyset$ and $\dot{x} \neq \emptyset$;
(e) $x \in \widehat{K} \iff x^+ \neq \emptyset, x^- \neq \emptyset$ and $\dot{x} = \emptyset$;
(f) $x \in \Delta \mathbb{R}^n \iff \dot{x} \neq \emptyset$.

**Lemma 4.1.** *Let $K$ be the positive cone of $\mathbb{R}^n$ with negative cone $-K$.*

(a) *For any $x \in \mathbb{R}^n$, $P_K(x)$ has the following representation*

$$P_K(x)_i = \begin{cases} x_i, & \text{if } i \in x^+, \\ 0, & \text{if } i \notin x^+, \end{cases} \text{ for } i = 1, 2, \ldots, n. \tag{4.1}$$

*In particular, we have*

(i) $P_K(x) = x$, *for any $x \in K$*;
(ii) $P_K(x) = \theta$, *for any $x \in -K$*;
(iii) $P_K(x) \in \partial K$, *for any $x \in \mathbb{R}^n \setminus K$*;
(iv) $P_K(x) \in \partial K \setminus \{\theta\}$, *for any $x \in \widehat{K}$*.

(b) *$P_K$ is positive homogeneous. For any $x \in \mathbb{R}^n$,*

$$P_K(\lambda x) = \lambda P_K(x), \text{ for any } \lambda \geq 0.$$

*Proof.* We write $P_K(x) = y = (y_1, y_2, \ldots, y_n) \in K$ with

$$y_i = \begin{cases} x_i, & \text{if } i \in x^+, \\ 0, & \text{if } i \notin x^+, \end{cases} \text{ for } i = 1, 2, \ldots, n.$$

For any $z = (z_1, z_2, \ldots, z_n) \in K$, we calculate

$$\langle x - y, y - z \rangle$$
$$= \langle (x_1, x_2, \ldots, x_n) - (y_1, y_2, \ldots, y_n), (y_1, y_2, \ldots, y_n) - (z_1, z_2, \ldots, z_n) \rangle$$
$$= \sum_{i=1}^n (x_i - y_i)(y_i - z_i)$$
$$= \sum_{i \in x^+} (x_i - y_i)(y_i - z_i) + \sum_{i \notin x^+} (x_i - y_i)(y_i - z_i)$$
$$= \sum_{i \in x^+} 0(y_i - z_i) + \sum_{i \notin x^+} x_i(-z_i)$$
$$= \sum_{i \notin x^+} x_i(-z_i) \geq 0.$$

By the basic variational principle of $P_K$, this proves $y = P_K(x)$. We see those parts (i), (ii) and (iii) follow from (4.1) immediately. We only show (iv). For any $x \in \widehat{K}$, by (4.1), we have that

$x^- \neq \emptyset$ and $x^+ \neq \emptyset$. By the representation of $P_K(x)$ in (4.1), it follows that $P_K(x) \in \partial K \setminus \{\theta\}$. □

For any given fixed $x \in \widehat{K}$, we define a mapping $b(x; \cdot): \mathbb{R}^n \to \mathbb{R}^n$, for any $w \in \mathbb{R}^n$, by

$$(b(x; w))_i = \begin{cases} w_i, & \text{if } i \in x^+, \\ 0, & \text{if } i \notin x^+, \end{cases} \text{ for } i = 1, 2, \ldots, n. \tag{4.2}$$

**Lemma 4.2**. *For any fixed $x \in \widehat{K}$, as defined in (4.1), $b(x; \cdot)$ satisfies*

(a) $b(x; \cdot): \mathbb{R}^n \to \mathbb{R}^n$ *is a linear and continuous mapping*;
(b) $b(x; w) \in \Delta \mathbb{R}^n$, *for any $w \in \mathbb{R}^n \setminus K$*.

*Proof*. The proof of this lemma is trivial and it is omitted here. □

For any given $x \in \Delta \mathbb{R}^n$, we define a mapping $d(x; \cdot): \mathbb{R}^n \to \mathbb{R}^n$, for any $w \in \mathbb{R}^n$, by

$$\begin{aligned}
(d(x; w))_i &= w_i, &&\text{for } i \in x^+, \\
(d(x; w))_i &= 0, &&\text{for } i \in x^-, \\
\text{and} \quad (d(x; w))_i &= \begin{cases} w_i, & \text{if } w_i > 0, \\ 0, & \text{if } w_i \leq 0, \end{cases} &&\text{for } i \in \dot{x}.
\end{aligned} \tag{4.3}$$

**Lemma 4.3**. *For any fixed $x \in \Delta \mathbb{R}^n$, as defined in (4.2), $d(x; \cdot)$ satisfies*

(a) $d(x; \cdot)$ *is a non-liner mapping from $\mathbb{R}^n$ to $\mathbb{R}^n$*;
(b) $d(\theta; w) = P_K(w)$, *for any $w \in \mathbb{R}^n$*.

Proof. The proof of this lemma is straight forward and it is omitted here. □

By using the results of Lemma 4.1, we study the strict Fréchet differentiability of the metric projection operator $P_K$.

**Theorem 4.4**. *Let $K$ be the positive cone of $\mathbb{R}^n$ with negative cone $-K$. Then, the metric projection operator $P_K$ has the following Fréchet differentiability properties.*

(i) $P_K$ *is strict Fréchet differentiable on $K^o$ satisfying $\nabla P_K(x) = I_{\mathbb{R}^n}$, for any $x \in K^o$, so*

$$\nabla P_K(x)(y) = y, \text{ for any } y \in \mathbb{R}^n;$$

(ii) $P_K$ *is strict Fréchet differentiable on $(-K)^o$ satisfying $\nabla P_K(x) = \theta$, for any $x \in (-K)^o$,*

$$\nabla P_K(x)(y) = \theta, \text{ for any } y \in \mathbb{R}^n;$$

(iii) $P_K$ *is strictly Fréchet differentiable on $\widehat{K}$ such that, for any $x \in \widehat{K}$, so*

$$\nabla P_K(x)(y) = b(x; y), \text{ for any } y \in \mathbb{R}^n;$$

(iv) *For subset $\Delta \mathbb{R}^n$, we have*

(a) $P_K$ *is Gâteaux directionally differentiable on $\Delta \mathbb{R}^n$ such that, for any $x \in \Delta \mathbb{R}^n$,*

$$P'_K(x)(w) = d(x; w), \quad \text{for any } w \in \mathbb{R}^n \setminus \{\theta\}.$$

(b) $P_K$ is not Fréchet differentiable at any point in $\Delta \mathbb{R}^n$, that is,

$$\nabla P_K(x) \text{ does not exist, for any } x \in \Delta \mathbb{R}^n.$$

*Proof.* Part (i) follows from Proposition 2.3 and Part (ii) follows from Proposition 2.2.

Proof of (a) in part (iii). For any $x \in \widehat{K}$, by the definition of Fréchet differentiability of the metric projection operator $P_K$, we consider the following limit.

$$\lim_{u \to x} \frac{P_K(u) - P_K(x) - b(x; u-x)}{\|u - x\|}. \tag{4.4}$$

For $x \in \widehat{K}$, by definition, $x^+ \neq \emptyset$, $x^- \neq \emptyset$ and $x^+ \cup x^- = \{1, 2, \ldots, n\}\}$, which implies that $\dot{x} = \emptyset$. Let $\delta_x = \frac{1}{4}\min\{|x_i|: i = 1, 2, \ldots, n\}$. Then $\delta_x > 0$ such that, for any $u, v \in \mathbb{R}^n$, if $\|u - x\| < \delta_x$ and $\|v - x\| < \delta_x$, we have

(a) $u^+ = x^+$ and $u^- = x^-$;
(b) $v^+ = x^+$ and $v^- = x^-$.

By Lemma 4.1, this implies that for any $u \in \mathbb{R}^n$, if $\|u - x\| < \delta_x$ and $\|v - x\| < \delta_x$, then, for $i = 1, 2, \ldots, n$, we have

$$P_K(u)_i = \begin{cases} u_i, & \text{if } i \in x^+, \\ 0, & \text{if } i \notin x^+; \end{cases} \quad \text{and} \quad P_K(v)_i = \begin{cases} v_i, & \text{if } i \in x^+, \\ 0, & \text{if } i \notin x^+. \end{cases} \tag{4.5}$$

By the definition of $b(x; \cdot)$, we have

$$b(x; u-v)_i = \begin{cases} u_i - v_i, & \text{if } i \in x^+, \\ 0, & \text{if } i \notin x^+. \end{cases} \tag{4.6}$$

Substituting (4.5) and (4.6) into (4.4), we obtain

$$\lim_{(u,v) \to (x,x)} \frac{P_K(u) - P_K(v) - b(x; u-v)}{\|u - v\|}$$

$$= \lim_{\substack{(u,v) \to (x,x) \\ \|u-x\| < \delta_x, \|v-x\| < \delta_x}} \frac{P_K(u) - P_K(v) - b(x; u-v)}{\|u - v\|}$$

$$= \lim_{\substack{(u,v) \to (x,x) \\ \|u-x\| < \delta_x, \|v-x\| < \delta_x}} \frac{\theta}{\|u - v\|}$$

$$= \theta.$$

This proves part (iii).

Proof of (iv). We prove (a) of part (iv). For any given $x \in \Delta \mathbb{R}^n$, we prove that $P_K$ is Gâteaux directionally differentiable at point $x$ such that

$$P'_K(x)(w) = d(x; w), \text{ for any } w \in \mathbb{R}^n \setminus \{\theta\}.$$

For this arbitrarily given $x \in \Delta \mathbb{R}^n$ with $\dot{x} \neq \emptyset$, for any fixed $w \in \mathbb{R}^n$ with $w \neq \theta$, we consider the following limit.

$$\lim_{t \downarrow 0} \frac{P_K(x+tw) - P_K(x)}{t}. \tag{4.7}$$

By (4.2) and (4.1), we have

$$d(\theta; w) = P_K(w), \text{ for any } w \in \mathbb{R}^n \setminus \{\theta\}.$$

This implies that, if $x = \theta$, then, by the property that $P_K(tw) = tP_K(w)$, for $t > 0$, we have

$$P_K(\theta + tw) - P_K(\theta) = tP_K(w) - \theta = td(\theta; w). \tag{4.8}$$

Substituting (4.8) into (4.7), we obtain

$$P'_K(\theta)(w) = \lim_{t \downarrow 0} \frac{P_K(\theta+tw) - P_K(\theta)}{t} = \lim_{t \downarrow 0} \frac{td(\theta; w)}{t} = d(\theta; w).$$

This proves (a) of (iv) for $x = \theta$. Next, we suppose that $x \neq \theta$. Let

$$\delta_x = \frac{1}{4}\min\{|x_i| : x_i \neq 0, i = 1, 2, \ldots, n\},$$

and

$$\Delta_w = \max\{|w_i| : i = 1, 2, \ldots, n\}, \text{ for } w \in \mathbb{R}^n \setminus \{\theta\}.$$

In the limit (4.7), if $0 < t < \frac{\delta_x}{2\Delta_w}$, then, for $i = 1, 2, \ldots, n$, we have

$$x_i > 0 \quad \Rightarrow \quad x_i + tw_i > 0,$$

and

$$x_i < 0 \quad \Rightarrow \quad x_i + tw_i < 0.$$

This implies that, if $0 < t < \frac{\delta_x}{2\Delta_w}$, then,

$$P_K(x + tw)_i = x_i + tw_i, \quad \text{for } i \in x^+$$

$$P_K(x + tw)_i = 0, \quad \text{for } i \in x^-,$$

and

$$P_K(x + tw)_i = \begin{cases} tw_i, & \text{if } w_i > 0, \\ 0, & \text{if } w_i \leq 0, \end{cases} \quad \text{for } i \in \dot{x}.$$

This implies that, for $i = 1, 2, \ldots, n$, we have

$$P_K(x + tw)_i - P_K(x)_i = x_i + tw_i - x_i = tw_i, \quad \text{for } i \in x^+,$$

$$P_K(x + tw)_i - P_K(x)_i = 0 - 0 = 0, \qquad \text{for } i \in x^-,$$

and

$$P_K(x + tw)_i - P_K(x)_i = \begin{cases} tw_i, & \text{if } w_i > 0, \\ 0, & \text{if } w_i \leq 0, \end{cases} \quad \text{for } i \in \dot{x}.$$

By the definition (4.2) of $d(x; \cdot)$, this implies

$$P_K(x + tw) - P_K(x) = d(x; tw) = td(x; w). \tag{4.9}$$

Substituting (4.9) into (4.7), we obtain

$$P_K'(x)(w) = \lim_{t \downarrow 0} \frac{P_K(x+tw) - P_K(x)}{t} = \lim_{t \downarrow 0} \frac{td(x; w)}{t} = d(x; w).$$

This proves (a) of part (iv).

Then, we prove (b) of part (iv). We prove that $P_K$ is not Fréchet differentiable at any point in $\Delta \mathbb{R}^n$. For an arbitrarily given $x \in \Delta \mathbb{R}^n$, there is $k \in \{1, 2, \ldots, n\}$ such that $x_k = 0$. Assume, by the way of contradiction, that $P_K$ is Fréchet differentiable at $x \in \Delta \mathbb{R}^n$. Then, there is a linear and continuous mapping $A(x): \mathbb{R}^n \to \mathbb{R}^n$ such that

$$\lim_{u \to x} \frac{P_K(u) - P_K(x) - A(x)(u-x)}{\|u-x\|} = \theta. \tag{4.10}$$

Let

$$z = (z_1, z_2, \ldots, z_n) \in \mathbb{R}^n \text{ with } z_i = 0, \text{ for } i \neq k \text{ and } z_k = 1.$$

In the limit (4.10), we take a special directional line segment for $u$ approaching to $x$ by $u = (u_1, u_2, \ldots, u_n)$ satisfying

$$u_i = x_i, \text{ for } i \neq k \quad \text{and} \quad u_k = -t, \text{ for } t > 0.$$

By the formula (4.1), we have

$$u - x = -tz \quad \text{and} \quad P_K(u) - P_K(x) = \theta. \tag{4.11}$$

Substituting (4.11) into (4.10), we have

$$\theta = \lim_{u \to x} \frac{P_K(u) - P_K(x) - A(x)(u-x)}{\|u-x\|}$$

$$= \lim_{t \downarrow 0} \frac{\theta - A(x)(-tz)}{\|-tz\|}$$

$$= A(x)(z).$$

This implies

$$A(x)(z) = \theta. \tag{4.12}$$

In the limit (4.10), we similarly take another special directional line segment for $v$ approaching to $x$ by $v = (v_1, v_2, \ldots, v_n)$ satisfying

$$v_i = x_i, \text{ for } i \neq k \quad \text{and} \quad v_k = t, \text{ for } t > 0.$$

By the formula (4.1), we calculate

$$v - x = tz \quad \text{and} \quad P_K(v) - P_K(x) = tz. \qquad (4.13)$$

Substituting (4.13) into (4.10), we have

$$\theta = \lim_{v \to x} \frac{P_K(v) - P_K(x) - A(x)(v-x)}{\|v-x\|}$$

$$= \lim_{t \downarrow 0} \frac{tz - A(x)(tz)}{\|-tw\|}$$

$$= z - A(x)(z).$$

This implies

$$A(x)(z) = z.$$

Since $z \neq \theta$, this contradicts to (4.12). This proves that $P_K$ is not Fréchet differentiable at this arbitrarily given point $x \in \Delta \mathbb{R}^n$, this proves (b) of part (iv). □

## 5. The Fréchet non-differentiability and Gâteaux directional differentiability of the metric projection onto the positive cone in $l_2$

In this section, we consider the real Hilbert space $l_2$ with norm $\|\cdot\|$, with inner product $\langle \cdot, \cdot \rangle$ and with the origin $\theta$. Let $\mathbb{N}$ denote the set of all positive integers. Let $N$ be a nonempty subset of $\mathbb{N}$ with complement $\bar{N}$. We define some subsets in $l_2$:

$$\mathbb{K} = \{x = (x_1, x_2, \ldots) \in l_2 : x_i \geq 0, \text{ for all } i \in \mathbb{N}\},$$

$$\mathbb{K}^+ = \{x = (x_1, x_2, \ldots) \in \mathbb{K} : x_i > 0, \text{ for all } i \in \mathbb{N}\},$$

$$\mathbb{K}^- = -\mathbb{K}^+ = \{x = (x_1, x_2, \ldots) \in \mathbb{K} : x_i < 0, \text{ for all } i \in \mathbb{N}\},$$

$$\widehat{\mathbb{K}} = \{x = (x_1, x_2, \ldots) \in l_2 : |x_i| > 0, \text{ for all } i \in \mathbb{N} \text{ and}$$
$$\text{there are at least one pair } j, k \in \mathbb{N} \text{ with } x_j x_k < 0\}.$$

$\mathbb{K}$ is a pointed closed and convex cone that is called the positive cone in $l_2$. In contrast with Euclidean spaces studied in the previous section, the interior of positive cone $\mathbb{K}$ in $l_2$ is empty. We prove it by the following lemma.

**Lemma 5.1.** *The interior of the positive cone $\mathbb{K}$ (the negative cone $\mathbb{K}^-$) in $l_2$ is empty.*

*Proof.* Let $x = (x_1, x_2, \ldots) \in \mathbb{K}$. For any $\varepsilon > 0$, there is a positive integer $m$ such that

$$\sum_{i=n}^{\infty} x_i^2 < \frac{1}{4}\varepsilon^2, \text{ for any } n > m.$$

Let $y = (y_1, y_2, \ldots) \in l_2$ satisfying

$$y_i = \begin{cases} x_i, & \text{if } i \le m, \\ -\frac{1}{2}\varepsilon, & \text{if } i = m+1, \\ 0, & \text{if } i > m+1. \end{cases}$$

One can check that $\|x - y\| < \varepsilon$ and $y \notin \mathbb{K}$. $\square$

Let $P_\mathbb{K}$ be the metric projection operator from $l_2$ onto $\mathbb{K}$. Similar, to Lemma 4.1, $P_\mathbb{K}$ has the following representation properties.

**Lemma 5.2.** *Let $\mathbb{K}$ be the positive cone of $l_2$. $P_\mathbb{K}$ has the following properties.*

(a) *For any $x \in l_2$, $P_\mathbb{K}(x)$ is represented as follows*

$$P_\mathbb{K}(x)_i = \begin{cases} x_i, & \text{if } x_i > 0, \\ 0, & \text{if } x_i \le 0, \end{cases} \text{ for } i = 1, 2, \ldots . \tag{5.1}$$

(b) *$P_\mathbb{K}$ is positive homogeneous. For any $x \in l_2$,*

$$P_\mathbb{K}(\lambda x) = \lambda P_\mathbb{K}(x), \text{ for any } \lambda \ge 0. \tag{5.2}$$

*Proof.* The proof of part (a) is similar to the proof of part (a) in Theorem 4.3, which is proved by using the basic variational principle of $P_\mathbb{K}$. Part (b) follows from part (a). $\square$

Similarly, to (4.2) in section 4, for any given fixed $x \in \widehat{\mathbb{K}}$, we define a mapping $B(x; \cdot): l_2 \to l_2$, for any $w \in l_2$, by

$$(B(x; w))_i = \begin{cases} w_i, & \text{if } x_i > 0, \\ 0, & \text{if } x_i \le 0, \end{cases} \text{ for all } i \in \mathbb{N}.$$

In contrast with Theorem 4.4 in the previous section, in the next theorem, we prove that $P_\mathbb{K}$ is not Fréchet differentiable on both $\mathbb{K}^+$ and $\mathbb{K}^-$.

**Theorem 5.3.** *Let $\mathbb{K}$ be the positive cone of $l_2$. Then, $P_\mathbb{K}$ has the following properties.*

(i) *In $\mathbb{K}^+$, we have*

(a) *$P_\mathbb{K}$ is not Fréchet differentiable at any point in $\mathbb{K}^+$. that is,*

$$\nabla P_\mathbb{K}(x) \text{ does not exist, for any } x \in \mathbb{K}^+.$$

(b) *$P_\mathbb{K}$ is Gâteaux directionally differentiable on $\mathbb{K}^+$ such that, for any $x \in \mathbb{K}^+$,*

$$P'_\mathbb{K}(x)(w) = w, \text{ for any } w \in l_2\setminus\{\theta\}.$$

(ii) *In $\mathbb{K}^-$, we have*

(a) *$P_\mathbb{K}$ is not Fréchet differentiable at any point in $\mathbb{K}^-$. that is,*

$$\nabla P_\mathbb{K}(x) \text{ does not exist, for any } x \in \mathbb{K}^-.$$

(b)     $P_\mathbb{K}$ is Gâteaux directionally differentiable on $\mathbb{K}^-$ such that, for any $x \in \mathbb{K}^-$,

$$P'_\mathbb{K}(x)(w) = \theta, \quad \text{for any } w \in \mathbb{R}^n \backslash \{\theta\}.$$

(iii)     In $\widehat{\mathbb{K}}$, we have

(a)     $P_\mathbb{K}$ is not Fréchet differentiable on $\widehat{\mathbb{K}}$, that is,

$$\nabla P_\mathbb{K}(x) \text{ does not exist, for any } x \in \widehat{\mathbb{K}}.$$

(b)     $P_\mathbb{K}$ is Gâteaux directionally differentiable on $\widehat{\mathbb{K}}$ such that, for any $x \in \widehat{\mathbb{K}}$,

$$P'_\mathbb{K}(x)(w) = B(x; w), \quad \text{for any } w \in l_2 \backslash \{\theta\}.$$

*Proof.* Since the proofs of (a)'s of (i), (ii) and (iii) are similar and the proofs of (b)'s of (i), (ii) and (iii) are similar, so, we first prove (a)'s in (i), (ii) and (iii). Then, we prove (b)'s.

Proof of (a) of part (i). For arbitrarily given $x \in \mathbb{K}^+$, assume, by the way of contradiction, that $P_\mathbb{K}$ is Fréchet differentiability at $x$. Then there is a linear and continuous mapping $A(x): l_2 \to l_2$, such that

$$\lim_{u \to x} \frac{P_\mathbb{K}(u) - P_\mathbb{K}(x) - A(x)(u-x)}{\|u - x\|} = \theta. \tag{5.3}$$

For any $n \in \mathbb{N}$, we write $w(n) \in l_2$ by

$$w(n)_i = \begin{cases} 0, & \text{if } i \neq n, \\ 1, & \text{if } i = n, \end{cases} \text{ for } i \in \mathbb{N}.$$

We define $u(n), v(n) \in l_2$ by

$$u(n)_i = \begin{cases} x_i, & \text{if } i \neq n, \\ -x_n, & \text{if } i = n, \end{cases} \text{ for } i \in \mathbb{N},$$

and

$$v(n)_i = \begin{cases} x_i, & \text{if } i \neq n, \\ -2x_n, & \text{if } i = n, \end{cases} \text{ for } i \in \mathbb{N}.$$

It is clear that

$$u(n) - x = -2x_n w(n) \to \theta, \text{ as } n \to \infty, \tag{5.4}$$

$$P_\mathbb{K}(u(n)) - P_\mathbb{K}(x) = -x_n w(n), \tag{5.5}$$

$$v(n) - x = -3x_n w(n) \to \theta, \text{ as } n \to \infty. \tag{5.7}$$

$$P_\mathbb{K}(v(n)) - P_\mathbb{K}(x) = -x_n w(n). \tag{5.8}$$

In the limit (5.3), we take a directional line segment $u = u(n)$, for $n \to \infty$. By (5.1), (5.2), (5.4) and (5.5), and by the linearity and continuity of the mapping $A(x)$, we have

$$\theta = \lim_{n\to\infty} \frac{P_{\mathbb{K}}(u(n))-P_{\mathbb{K}}(x)-A(x)(u(n)-x)}{\|u(n)-x\|}$$

$$= \lim_{n\to\infty} \frac{-x_n w(n) - A(x)(-2x_n w(n))}{\|-2x_n w(n)\|}$$

$$= \lim_{n\to\infty} \frac{-x_n w(n) + 2x_n A(x)(w(n))}{2x_n}$$

$$= \lim_{n\to\infty} \frac{-w(n) + 2A(x)(w(n))}{2}.$$

This implies

$$\lim_{n\to\infty}\left(-w(n) + 2A(x)(w(n))\right) = \theta. \tag{5.9}$$

In the limit (5.3), we take another directional line segment $v = v(n)$, for $n \to \infty$. By (5.1), (5.2), (5.7) and (5.8), and by the linearity and continuity of the mapping $A(x)$, we have

$$\theta = \lim_{n\to\infty} \frac{P_{\mathbb{K}}(v(n))-P_{\mathbb{K}}(x)-A(x)(v(n)-x)}{\|v(n)-x\|}$$

$$= \lim_{n\to\infty} \frac{-x_n w(n) - A(x)(-3x_n w(n))}{\|-3x_n w(n)\|}$$

$$= \lim_{n\to\infty} \frac{-x_n w(n) + 3x_n A(x)(w(n))}{3x_n}$$

$$= \lim_{n\to\infty} \frac{-w(n) + 3A(x)(w(n))}{3}.$$

This implies

$$\lim_{n\to\infty}\left(-w(n) + 3A(x)(w(n))\right) = \theta. \tag{5.10}$$

2×(5.10) minus 3×(5.9), we get

$$\lim_{n\to\infty} w(n) = \theta. \tag{5.11}$$

(5.11) is a contradiction, which proves (a) of part (i).

Proof of (a) in (ii). The proof of (a) in part (ii) is similar to the proof of (a) in part (i).

For arbitrarily given $y \in \mathbb{K}^-$, assume, by the way of contradiction, that $P_{\mathbb{K}}$ is Fréchet differentiability at $y$. Then there is a linear and continuous mapping $B(y): l_2 \to l_2$, such that

$$\lim_{u\to x} \frac{P_{\mathbb{K}}(u)-P_{\mathbb{K}}(y)-B(y)(u-y)}{\|u-y\|} = \theta. \tag{5.3}$$

Similarly, to the proof of (i), for any $n \in \mathbb{N}$, we write $w(n) \in l_2$ by

$$w(n)_i = \begin{cases} 0, & \text{if } i \neq n, \\ 1, & \text{if } i = n \end{cases}, \text{ for } i \in \mathbb{N}.$$

We define $u(n)$, $v(n) \in l_2$ by

$$u(n)_i = \begin{cases} y_i, & \text{if } i \neq n, \\ -y_n, & \text{if } i = n \end{cases}, \text{ for } i \in \mathbb{N},$$

and

$$v(n)_i = \begin{cases} y_i, & \text{if } i \neq n, \\ -2y_n, & \text{if } i = n \end{cases}, \text{ for } i \in \mathbb{N}.$$

It is clear that

$$u(n) - y = -2y_n w(n) \to \theta, \text{ as } n \to \infty, \tag{5.12}$$

$$P_\mathbb{K}(u(n)) - P_\mathbb{K}(y) = -y_n w(n), \tag{5.13}$$

$$v(n) - y = -3y_n w(n) \to \theta, \text{ as } n \to \infty. \tag{5.14}$$

$$P_\mathbb{K}(v(n)) - P_\mathbb{K}(y) = -2y_n w(n). \tag{5.15}$$

In the limit (5.3), we take a sequential approaching $u = u(n)$, for $n \to \infty$. By (5.1), (5.2), (5.12) and (5.13), and by the linearity and continuity of the mapping $B(y)$, we have

$$\theta = \lim_{n \to \infty} \frac{P_\mathbb{K}(u(n)) - P_\mathbb{K}(y) - B(y)(u(n) - y)}{\|u(n) - y\|}$$

$$= \lim_{n \to \infty} \frac{-y_n w(n) - B(y)(-2y_n w(n))}{\|-2y_n w(n)\|}$$

$$= \lim_{n \to \infty} \frac{-y_n w(n) + 2y_n B(y)(w(n))}{-2y_n}$$

$$= \lim_{n \to \infty} \frac{w(n) - 2B(y)(w(n))}{2}.$$

This implies

$$\lim_{n \to \infty} \left( w(n) - 2B(y)(w(n)) \right) = \theta. \tag{5.16}$$

In the limit (5.3), we take a sequence approaching $v = v(n)$, for $n \to \infty$. By (5.1), (5.2), (5.14) and (5.15), and by the linearity and continuity of the mapping $B(y)$, we have

$$\theta = \lim_{n \to \infty} \frac{P_\mathbb{K}(v(n)) - P_\mathbb{K}(y) - B(y)(v(n) - y)}{\|u(n) - y\|}$$

$$= \lim_{n \to \infty} \frac{-2y_n w(n) - B(y)(-3y_n w(n))}{\|-3y_n w(n)\|}$$

$$= \lim_{n \to \infty} \frac{-2y_n w(n) + 3y_n B(y)(w(n))}{-3y_n}$$

$$= \lim_{n\to\infty} \frac{2w(n) - 3B(y)(w(n))}{3}.$$

This implies

$$\lim_{n\to\infty} \left(2w(n) - 3B(y)(w(n))\right) = \theta. \tag{5.17}$$

2×(5.17) minus 3×(5.16), we get

$$\lim_{n\to\infty} w(n) = \theta. \tag{5.11}$$

(5.11) is a contradiction, which proves (a) of part (ii).

Proof (a) in part (iii). For any given $x \in \widehat{\mathbb{K}}$, by the definition of $\widehat{\mathbb{K}}$, we have

$$\{i \in \mathbb{N}: x_i > 0\} \cup \{i \in \mathbb{N}: x_i < 0\} = \mathbb{N}.$$

This implies that, at least one of $\{i \in \mathbb{N}: x_i > 0\}$ and $\{i \in \mathbb{N}: x_i < 0\}$ is infinite. If $\{i \in \mathbb{N}: x_i > 0\}$ is infinite, then, then proof of (a) in part (iii) is very similar to the proof of (a) in part (i). If $\{i \in \mathbb{N}: x_i < 0\}$ is infinite, then, then proof of (a) in part (iii) is very similar to the proof of (a) in part (ii). So, the proof of (a) in part (iii) is omitted here.

Next, we prove (b) in part (i). For an arbitrarily given $x \in \mathbb{K}^+$, we prove

$$P'_{\mathbb{K}}(x)(w) = w, \text{ for any } w \in l_2\backslash\{\theta\}.$$

This is equivalent to prove

$$\lim_{t\downarrow 0} \frac{P_{\mathbb{K}}(x+tw) - P_{\mathbb{K}}(x)}{t} = w, \text{ for any } w \in l_2\backslash\{\theta\}.$$

For any fixed $w \in l_2\backslash\{\theta\}$, we have $\|w\| > 0$. It is clear to see that $\|w\| \geq |w_i|$, for all $i \in \mathbb{N}$. For any given $\varepsilon > 0$, for this arbitrarily given $w \in l_2\backslash\{\theta\}$, there is a positive integer $n$ such that

$$(\textstyle\sum_{i>n} w_i^2)^{\frac{1}{2}} < \frac{\varepsilon}{2}.$$

For this given $n$, let $\delta = \min\{x_i: i = 1, 2, \ldots, n\}$. Then, if $0 < t < \frac{\delta}{2\|w\|}$, then

$$x_i + tw_i > 0, \text{ for all } i = 1, 2, \ldots, n.$$

This implies that, if $0 < t < \frac{\delta}{2\|w\|}$, then

$$P_{\mathbb{K}}(x+tw)_i = x_i + tw_i, \text{ for all } i = 1, 2, \ldots, n.$$

Then, for $0 < t < \frac{\delta}{2\|w\|}$, we estimate

$$\left\|\frac{P_{\mathbb{K}}(x+tw) - P_{\mathbb{K}}(x)}{t} - w\right\|^2$$

$$= \frac{\|P_{\mathbb{K}}(x+tw)-P_{\mathbb{K}}(x)-tw\|^2}{t^2}$$

$$= \frac{\sum_{i=1}^{\infty}(P_{\mathbb{K}}(x+tw)-P_{\mathbb{K}}(x)-tw)_i^2}{t^2}$$

$$= \frac{\sum_{i=1}^{n}(P_{\mathbb{K}}(x+tw)-P_{\mathbb{K}}(x)-tw)_i^2 + \sum_{i=n+1}^{\infty}(P_{\mathbb{K}}(x+tw)-P_{\mathbb{K}}(x)-tw)_i^2}{t^2}$$

$$= \frac{\sum_{i=1}^{n} 0 + \sum_{i=n+1}^{\infty}(P_{\mathbb{K}}(x+tw)-P_{\mathbb{K}}(x)-tw)_i^2}{t^2}$$

$$= \frac{\sum_{i>n, x_i+tw_i>0}(P_{\mathbb{K}}(x+tw)-P_{\mathbb{K}}(x)-tw)_i^2 + \sum_{i>n, x_i+tw_i \leq 0}(P_{\mathbb{K}}(x+tw)-P_{\mathbb{K}}(x)-tw)_i^2}{t^2}$$

$$= \frac{\sum_{i>n, x_i+tw_i>0} 0 + \sum_{i>n, x_i+tw_i \leq 0}(P_{\mathbb{K}}(x+tw)-P_{\mathbb{K}}(x)-tw)_i^2}{t^2}$$

$$= \frac{\sum_{i>n, x_i+tw_i \leq 0}(0-P_{\mathbb{K}}(x)-tw)_i^2}{t^2}$$

$$= \frac{\sum_{i>n, x_i+tw_i \leq 0}(P_{\mathbb{K}}(x)+tw)_i^2}{t^2}$$

$$= \frac{\sum_{i>n, x_i+tw_i \leq 0}(x+tw)_i^2}{t^2}$$

$$\leq \frac{2\sum_{i>n, x_i+tw_i \leq 0}(x_i^2+t^2 w_i^2)}{t^2}$$

$$\leq \frac{2\sum_{i>n, x_i+tw_i \leq 0}(t^2 w_i^2+t^2 w_i^2)}{t^2}$$

$$\leq 4 \sum_{i>n} w_i^2$$

$$< \varepsilon^2.$$

This implies that, if $0 < t < \frac{\delta}{2\|w\|}$, then

$$\left\|\frac{P_{\mathbb{K}}(x+tw)-P_{\mathbb{K}}(x)}{t} - w\right\| < \varepsilon.$$

This proves (b) in part (i) that

$$\lim_{t \downarrow 0} \frac{P_{\mathbb{K}}(x+tw)-P_{\mathbb{K}}(x)}{t} = w, \text{ for } w \in l_2 \setminus \{\theta\}.$$

Next, we prove (b) in part (ii). For an arbitrarily given $x \in \mathbb{K}^-$, we prove

$$P'_{\mathbb{K}}(x)(w) = \theta, \text{ for any } w \in l_2 \setminus \{\theta\}.$$

This is equivalent to prove

$$\lim_{t \downarrow 0} \frac{P_\mathbb{K}(x+tw) - P_\mathbb{K}(x)}{t} = \theta, \text{ for any } w \in l_2 \setminus \{\theta\}.$$

For any fixed $w \in l_2 \setminus \{\theta\}$, we have $\|w\| > 0$. It is clear to see that $\|w\| \geq |w_i|$, for all $i \in \mathbb{N}$. For any given $\varepsilon > 0$, for this arbitrarily given $w \in l_2 \setminus \{\theta\}$, there is a positive integer $n$ such that

$$(\Sigma_{i>n} w_i^2)^{\frac{1}{2}} < \frac{\varepsilon}{2}.$$

For this given $n$, let $\delta = \min\{-x_i : i = 1, 2, \ldots, n\}$. Then, if $0 < t < \frac{\delta}{2\|w\|}$, then

$$x_i + tw_i < 0, \text{ for all } i = 1, 2, \ldots, n.$$

This implies that, if $0 < t < \frac{\delta}{2\|w\|}$, then

$$P_\mathbb{K}(x+tw)_i = 0, \text{ for all } i = 1, 2, \ldots, n.$$

Then, for $0 < t < \frac{\delta}{2\|w\|}$, we estimate

$$\left\| \frac{P_\mathbb{K}(x+tw) - P_\mathbb{K}(x)}{t} - \theta \right\|^2$$

$$= \frac{\|P_\mathbb{K}(x+tw) - P_\mathbb{K}(x)\|^2}{t^2}$$

$$= \frac{\|P_\mathbb{K}(x+tw)\|^2}{t^2}$$

$$= \frac{\Sigma_{i=1}^\infty (P_\mathbb{K}(x+tw))_i^2}{t^2}$$

$$= \frac{\Sigma_{i=1}^n (P_\mathbb{K}(x+tw))_i^2 + \Sigma_{i=n+1}^\infty (P_\mathbb{K}(x+tw))_i^2}{t^2}$$

$$= \frac{\Sigma_{i=1}^n 0 + \Sigma_{i=n+1}^\infty (P_\mathbb{K}(x+tw))_i^2}{t^2}$$

$$= \frac{\Sigma_{i>n, x_i+tw_i>0} (P_\mathbb{K}(x+tw))_i^2 + \Sigma_{i>n, x_i+tw_i \leq 0} (P_\mathbb{K}(x+tw))_i^2}{t^2}$$

$$= \frac{\Sigma_{i>n, x_i+tw_i>0} (x+tw)_i^2 + \Sigma_{i>n, x_i+tw_i \leq 0} 0}{t^2}$$

$$= \frac{\Sigma_{i>n, x_i+tw_i>0} (x+tw)_i^2}{t^2}$$

$$\leq \frac{2 \Sigma_{i>n, x_i+tw_i>0} (x_i^2 + t^2 w_i^2)}{t^2}$$

$$\leq \frac{2 \Sigma_{i>n, x_i+tw_i>0} (t^2 w_i^2 + t^2 w_i^2)}{t^2}$$

$$\leq 4\sum_{i>n} w_i^2$$

$$< \varepsilon^2.$$

This implies that, if $0 < t < \frac{\delta}{2\|w\|}$, then

$$\left\|\frac{P_{\mathbb{K}}(x+tw)-P_{\mathbb{K}}(x)}{t}\right\| < \varepsilon.$$

This proves (b) in part (ii) that

$$\lim_{t\downarrow 0}\frac{P_{\mathbb{K}}(x+tw)-P_{\mathbb{K}}(x)}{t} = \theta, \text{ for } w \in l_2\backslash\{\theta\}.$$

Next, we prove (b) in part (iii). For an arbitrarily given $x \in \widehat{\mathbb{K}}$, we prove

$$P'_{\mathbb{K}}(x)(w) = B(x; w), \quad \text{for any } w \in l_2\backslash\{\theta\}.$$

This is equivalent to prove

$$\lim_{t\downarrow 0}\frac{P_{\mathbb{K}}(x+tw)-P_{\mathbb{K}}(x)}{t} = B(x; w), \text{ for any } w \in l_2\backslash\{\theta\}.$$

For any fixed $w \in l_2\backslash\{\theta\}$, we have $\|w\| > 0$. It is clear to see that $\|w\| \geq |w_i|$, for all $i \in \mathbb{N}$. For any given $\varepsilon > 0$, for this arbitrarily given $w \in l_2\backslash\{\theta\}$, there is a positive integer $n$ such that

$$(\sum_{i>n} w_i^2)^{\frac{1}{2}} < \frac{\varepsilon}{2}.$$

For this given $n$, let $\delta = \min\{|x_i|: i = 1, 2, \ldots, n\}$. Then, if $0 < t < \frac{\delta}{2\|w\|}$, then, for $i = 1, 2, \ldots, n$

$$x_i + tw_i > 0, \text{ if } x_i > 0,$$

and $$x_i + tw_i < 0, \text{ if } x_i < 0.$$

This implies that, if $0 < t < \frac{\delta}{2\|w\|}$, then, for $i = 1, 2, \ldots, n$, we have

$$P_{\mathbb{K}}(x+tw)_i = x_i + tw_i, \text{ if } x_i > 0,$$

$$P_{\mathbb{K}}(x+tw)_i = 0, \text{ if } x_i < 0.$$

Then, for $0 < t < \frac{\delta}{2\|w\|}$, we estimate

$$\left\|\frac{P_{\mathbb{K}}(x+tw)-P_{\mathbb{K}}(x)}{t} - B(x; w)\right\|^2$$

$$= \frac{\|P_{\mathbb{K}}(x+tw)-P_{\mathbb{K}}(x)-tB(x;w)\|^2}{t^2}$$

$$= \frac{\sum_{i=1}^{\infty}(P_{\mathbb{K}}(x+tw)-P_{\mathbb{K}}(x)-tB(x;w))_i^2}{t^2}$$

$$= \frac{\sum_{i=1}^{n}(P_{\mathbb{K}}(x+tw)-P_{\mathbb{K}}(x)-tB(x;w))_i^2 + \sum_{i=n+1}^{\infty}(P_{\mathbb{K}}(x+tw)-P_{\mathbb{K}}(x)-tB(x;w))_i^2}{t^2}$$

$$= \frac{\sum_{i=1}^{n}0 + \sum_{i=n+1}^{\infty}(P_{\mathbb{K}}(x+tw)-P_{\mathbb{K}}(x)-tB(x;w))_i^2}{t^2}$$

$$= \frac{\sum_{i>n, x_i+tw_i>0}(P_{\mathbb{K}}(x+tw)-P_{\mathbb{K}}(x)-tB(x;w))_i^2 + \sum_{i>n, x_i+tw_i\le 0}(P_{\mathbb{K}}(x+tw)-P_{\mathbb{K}}(x)-tB(x;w))_i^2}{t^2}$$

$$= \frac{\sum_{i>n, x_i+tw_i>0,\ x_i>0}(P_{\mathbb{K}}(x+tw)-P_{\mathbb{K}}(x)-tB(x;w))_i^2 + \sum_{i>n, x_i+tw_i>0,\ x_i\le 0}(P_{\mathbb{K}}(x+tw)-P_{\mathbb{K}}(x)-tB(x;w))_i^2}{t^2}$$

$$+ \frac{\sum_{i>n, x_i+tw_i\le 0,\ x_i>0}(P_{\mathbb{K}}(x+tw)-P_{\mathbb{K}}(x)-tB(x;w))_i^2 + \sum_{i>n, x_i+tw_i\le 0,\ x_i<0}(P_{\mathbb{K}}(x+tw)-P_{\mathbb{K}}(x)-tB(x;w))_i^2}{t^2}$$

$$= \frac{\sum_{i>n, x_i+tw_i>0,\ x_i>0}(x_i+tw_i-x_i-tw_i)^2 + \sum_{i>n, x_i+tw_i>0,\ x_i\le 0}(x_i+tw_i-0-t0)^2}{t^2}$$

$$+ \frac{\sum_{i>n, x_i+tw_i\le 0,\ x_i>0}(0-x_i-tw_i)^2 + \sum_{i>n, x_i+tw_i\le 0,\ x_i<0}(0-0-t0)^2}{t^2}$$

$$= \frac{\sum_{i>n, x_i+tw_i>0,\ x_i>0}0 + \sum_{i>n, x_i+tw_i>0,\ x_i\le 0}(x_i+tw_i)^2}{t^2} + \frac{\sum_{i>n, x_i+tw_i\le 0,\ x_i>0}(-x_i-tw_i)^2 + \sum_{i>n, x_i+tw_i\le 0,\ x_i<0}0}{t^2}$$

$$= \frac{\sum_{i>n, x_i+tw_i>0,\ x_i\le 0}(x_i+tw_i)^2 + \sum_{i>n, x_i+tw_i\le 0,\ x_i>0}(x_i+tw_i)^2}{t^2}$$

$$= \frac{2\sum_{i>n, x_i+tw_i>0,\ x_i\le 0}(x_i^2+t^2w_i^2) + 2\sum_{i>n, x_i+tw_i\le 0,\ x_i>0}(x_i^2+t^2w_i^2)}{t^2}$$

$$= \frac{2\sum_{i>n, x_i+tw_i>0,\ x_i\le 0}(t^2w_i^2+t^2w_i^2) + 2\sum_{i>n, x_i+tw_i\le 0,\ x_i>0}(t^2w_i^2+t^2w_i^2)}{t^2}$$

$$\le \frac{2\sum_{i>n}(t^2w_i^2+t^2w_i^2)}{t^2}$$

$$= 4\sum_{i>n}w_i^2$$

$$< \varepsilon^2.$$

This implies that, if $0 < t < \frac{\delta}{2\|w\|}$, then

$$\left\|\frac{P_{\mathbb{K}}(x+tw)-P_{\mathbb{K}}(x)}{t} - B(x;w)\right\| < \varepsilon.$$

This proves that

$$\lim_{t\downarrow 0}\frac{P_{\mathbb{K}}(x+tw)-P_{\mathbb{K}}(x)}{t} = B(x;w), \text{ for } w \in l_2\setminus\{\theta\}. \qquad \square$$

## 6. Conclusion and remarks

In section 3 of this paper, we study the strict Fréchet differentiability of the metric projection operator $P_{r\mathbb{B}}$ onto closed balls $r\mathbb{B}$ centered at the origin in Hilbert spaces. In [6, 7], the directional differentiability of the metric projection onto closed balls is studied in uniformly convex and uniformly smooth Banach spaces and Hilbert spaces, in which the considered balls have center at arbitrarily given point $c$ in the spaces. We believe that the results about the strict Fréchet differentiability of $P_{r\mathbb{B}}$ proved in Theorem 3.3 in this paper can be extended to metric projection $P_{\mathbb{B}(c,r)}$, which is onto closed balls with center at an arbitrarily given point $c$ in the spaces.

## Acknowledgments

The author is very grateful to Professor Boris Mordukhovich and Professor Simeon Reich for their kind communications, valuable suggestions and enthusiasm encouragements in the development stage of this paper.